\documentclass[11pt, eqsection]{article}
\usepackage{latexsym, graphicx}
\usepackage{amsmath}
\usepackage{amsfonts}
\usepackage{amssymb}
\usepackage{amssymb,amsmath,amsfonts,mathrsfs,amsxtra,amsthm}
\usepackage{graphicx,color,subfig}
\usepackage{accents,xspace}
\usepackage{dsfont}
\usepackage{amscd}
\newtheorem{Theorem}{Theorem}[section]
\newtheorem{Definition}[Theorem]{Definition}
\newtheorem{Proposition}[Theorem]{Proposition}
\newtheorem{Lemma}[Theorem]{Lemma}
\newtheorem{Corollary}[Theorem]{Corollary}

\newtheorem{Remark}[Theorem]{Remark}

\newenvironment{Proof}{\removelastskip\par\medskip
\noindent{\em Proof.} \rm}{\penalty-20\null\hfill$\square$\par\medbreak}

\def\N{\mathbb N }

\usepackage{amsfonts,layout}

\begin{document}

\pagenumbering{arabic}

\renewcommand{\theequation}{\thesection.\arabic{equation}}
\renewcommand{\theequation}{\thesubsection.\arabic{equation}}
\newcommand{\snl}{\ \\}
\newcommand{\pagref}[1]{page~\pageref{#1}}
\newcommand{\appref}[1]{appendix~\ref{#1}}
\newcommand{\chapref}[1]{chapter~\ref{#1}}
\newcommand{\tabref}[1]{table~\ref{#1}}
\newcommand{\figref}[1]{figure~\ref{#1}}
\newcommand{\secref}[1]{section~\ref{#1}}
\newcommand{\subsecref}[1]{subsection~\ref{#1}}
\newcommand{\ltfigure}[2]{\caption{#2} \label{#1}}

\newcommand{\m}[1]{{\(#1\)}}
\newcommand{\tilda}{\m{\sim}}
\newcommand{\mtilda}{\m{\sim}}
\newcommand{\msub}[2]{{\({#1}_{#2}\)}}
\newcommand{\msuper}[2]{{\({#1}^{#2}\)}}
\newcommand{\ucn}[1]{Ultracomputer Note \#{#1}}
\newcommand{\quotetitle}[1]{{\em ``#1''}}
\newcommand{\mla}{{$<=\ $}}
\newcommand{\smla}{{$<=$}}
\newcommand{\cons}{{$ \mid\ $}}

\newcommand{\samel}[1]{\mbox{#1}}
\newcommand{\ems}[1]{\samel{\em #1}}
\newcommand{\defeq}{\stackrel{\text{d̩f}}{=}}

\newcommand{\mypar}[1]{\paragraph{#1}\ \newline}

\newcommand{\rubsp}{\hspace{0em plus .25em}}
\newcommand{\rubspa}{\hspace{0em plus .1em}}

\def\squareforqed{\rule{1.5ex}{1.5ex}}

\def\qed{%
        \ifmmode\squareforqed\else{\unskip\nobreak\hfil%
        \penalty50\hskip1em\null\nobreak\hfil\squareforqed%
        \parfillskip=0pt\finalhyphendemerits=0\endgraf}\fi%
}

\let\finpreuve=\qed


\newcommand{\Nu}{\Upsilon}

\bibliographystyle{plain}

\title{{\bf \protect  Insensitizing exact controls for the scalar wave equation and exact controllability of  $2$-coupled cascade systems of PDE's by a single control. \\
}}
\author{{F}atiha {ALABAU-BOUSSOUIRA} \\ \\
        {\small {{U}niversit\'e de Lorraine, IECL UMR-CNRS 7502 and INRIA Team-project CORIDA}}\\
        {\small 
        57045 METZ {C}edex 1 ({F}rance)}}

\bibliographystyle{plain}
\date{}
\maketitle {}


\paragraph{AMS(MOS) subject classifications:} 34G10, 35B35, 35B37, 35L90, 93D15, 93D20.


\begin{abstract}
We study  the exact controllability, by a reduced number of controls, of coupled cascade systems of PDE's  and the existence of exact insensitizing controls for the scalar wave equation. 
We give a necessary and sufficient condition for the observability of abstract coupled cascade hyperbolic systems by a single observation, the  observation operator being either bounded or unbounded. Our proof extends the two-level energy method introduced in~\cite{sicon03, alaleau11} for symmetric coupled systems,  to cascade systems which are examples of non symmetric coupled systems.
In particular, we prove the observability of two coupled wave equations in cascade if the observation and coupling regions both satisfy the Geometric Control Condition (GCC) of Bardos Lebeau and Rauch~\cite{blr92}. 
By duality, this solves the exact controllability, by a single control, of $2$-coupled abstract cascade hyperbolic systems.
Using transmutation, we give null-controllability results for the multidimensional heat and Schr\"odinger $2$-coupled cascade systems under (GCC) and for any positive time. By our method, we can treat cases where the control and coupling coefficients have disjoint supports, 
partially solving an open question raised by de Teresa~\cite{DeT00}.  
Moreover we answer the question of the existence of  exact insensitizing locally distributed as well as boundary controls of scalar multidimensional wave equations, raised by J.-L. Lions~\cite{lions89} and later on by D\'ager~\cite{Dager06} and Tebou~\cite{tebou08}.
\end{abstract}

\section{Introduction}
The questions of the existence of insensitizing controls for scalar heat and wave equations are challenging issues. J.-L. Lions~\cite{lions89}  introduced this notion in 1989 to define controls which are robust to  small unknown perturbations on the initial data. It is by now well-known~\cite{lions89, Dager06, bodart-fabre95, DeT00, BGBPG04CPDE, BGBPG04SICON, DeTZ, tebou08, tebou2011} that the existence of such controls is equivalent to an exact controllability result for an associated $2$-coupled cascade of heat (resp. wave, \ldots) system, for which only one equation is controlled. Similar questions are considered for the Stokes equations in~\cite{guerrero}  and
for the Navier-Stokes equations in~\cite{gueyecras, gueyeNS}.  Hence the design of such insensitizing controls is related to the controllability, by a reduced number of controls, of coupled systems. Furthermore, many applicative issues in mechanics, biology or medecine lead also to similar controllability issues for coupled systems, which may be of parabolic, hyperbolic or mixed type. The coupling may also be more general than the cascade form. An increasing number of papers deals with these questions since then.

Coupled parabolic or diffusive control systems of order $2$ have the general form
\begin{equation}\label{parabolic}
\begin{cases}
e^{i \theta}y_t - \Delta y +\mathcal{C}y= Bv \,,\mbox{in }Q_T=\Omega\times (0,T)\,,\\
y=0 \,,\mbox{on } \Sigma_T=\partial \Omega \times (0,T)\,,\\
y(0,\cdot)=y_0(\cdot) \,,\mbox{in } \Omega \,,
\end{cases}
\end{equation}
with $\theta=0$ (resp. $\theta=\pi/2$) in the parabolic case  (resp. for Schr\"odinger case). Here $\Omega$ is an open non-empty subset in $\mathbb{R}^d$ with a smooth boundary 
$\Gamma$, $Y=(y_1, y_2)$ is the state to be controlled, $\mathcal{C}$
is a coupling bounded operator on $(L^2(\Omega))^2$, $B$ is either a bounded or unbounded control operator  acting on a single component of the above system, and $v$ is the scalar control . One can also consider the corresponding hyperbolic systems, obtained by replacing $e^{i \theta}y_t$ by $y_{tt}$ (together with appropriate initial conditions).

The above systems have received a lot of attention in the case of
cascade $2$-coupled parabolic systems. An example of such cascade system raises when $C$

\begin{equation}\label{cascade}
\mathcal{C}=
\begin{pmatrix}
0 &  \mathds{1}_O\\
0 & 0
\end{pmatrix}\,,
\end{equation}
with for instance $Bv=(0, v\mathds{1}_{\omega})^t$. Here $O$ and $\omega$ are open non empty subsets of $\Omega$ standing respectively for the coupling and control regions and $\mathds{1}_O$, the coupling coefficient, stands for the characteristic function of the set $O$. The coupling coefficient may be more generally a nonnegative (or nonpositive) function with a support that is included in $\Omega$. One may also consider other types of coupling operators such as, for instance, symmetric coupling operators 
\begin{equation}\label{symmcoup}
\mathcal{C}=
\begin{pmatrix}
0 &  c\\
c & 0
\end{pmatrix}\,,
\end{equation}
\noindent where $c\ge 0$ on $\Omega$.  One can also consider more general coupling operators.
\vskip 2mm

{\bf Some overview on the literature for controlled coupled systems} 

\vskip 2mm

Let us first present some of the results for parabolic equations. De Teresa~\cite{DeT00} considered $2$-coupled cascade parabolic systems with a single locally distributed control. Under the assumption that $\omega \cap O\neq \emptyset$, she
proved null controllability by a single control, and as a byproduct she obtained the existence of insensitizing controls for the scalar heat equation. In the case $\omega \cap O= \emptyset$, Kavian and de Teresa~\cite{DeTK10} proved a unique continuation result for a $2$-coupled cascade systems of parabolic equations. De Teresa and Zuazua~\cite{DeTZ} give further results, concerning the determination of the initial data for which insensitizing controls of the heat equation can be built. Null controllability results by a reduced number of controls, for more general coupled parabolic systems are proved mainly by means of Carleman estimates, and in two types of situations. The first one is devoted to constant (or time-dependent) coupling operators. The second one assumes locally distributed couplings, and locally distributed controls, but then with a non-empty intersection between the coupling and the control regions. We refer to~\cite{DeT00, AKBD06, AKBDGB09, FCGBDeT10, GBdT10, leautaud, CGR10, olive, mauffrey} and to the survey paper~\cite{AKBGBT11} for such results and to the references therein. Let us further mention an interesting result by Coron, Guerrero and Rosier~\cite{CGR10}, which proves local null controllability results for nonlinearly coupled $2$-systems of parabolic equations with a {\em nonlinear coupling} term. Such systems arise from applications to the control of chemical reaction-diffusion models. These nonlinear control results are based on the Coron's return method~\cite{coronbk07}. 

Hence it is a challenging issue to prove positive controllability results for coupled systems with a reduced number of controls, especially in situations for which the coupling and control regions do not meet. 
It seems that up to now the direct methods for coupled parabolic systems, mainly based on Carleman estimates, do not allow to deal with a coupling region that does not meet the control region, and this for boundary as well as locally distributed control problems. We shall see below that indirect methods, based for instance on controllability results for the corresponding hyperbolic system, may answer partially this issue.

\vskip 2mm

Let us now present some results of the literature for coupled hyperbolic systems. We shall first consider symmetric coupled systems
of wave equations (see~\cite{alacras01, sicon03}), that is

\begin{equation}\label{hyperbolic}
\begin{cases}
y_{1,tt} -\Delta y_1 + Cy_2= Bv \,,\mbox{in }Q_T=\Omega\times (0,T)\,,\\
y_{2,tt} -\Delta y_2 + C^{\ast}y_1= 0 \,,\mbox{in }Q_T=\Omega\times (0,T)\,,\\
y_i=0 \ i=1,2\,,\mbox{on } \Sigma_T=\partial \Omega \times (0,T)\,,\\
y_i(0,\cdot)=y_i^0(\cdot) \ i=1,2 \,,\mbox{in } \Omega \,.
\end{cases}
\end{equation}
\noindent Here only the first equation is controlled. One can also consider the corresponding abstract system in which the Dirichlet laplacian is replaced by a general self-adjoint coercive operator in a Hilbert space. In~\cite{sicon03}, we introduce an original method --named {\em two-level energy method}-- to prove positive controllability results. This method is set in a general abstract setting and relies on several properties. We use the property of conservation of the total energy of the adjoint system, and a time-independent observability inequality for a scalar equation with a source term.  Our method also relies on the idea to work in a weakened energy space for the unobserved component and to use a balance of energies between the unobserved and the observed components. Thanks to this, we prove observability and controllability results for coupled systems in the case of coercive bounded coupling operators $C$ (case of globally distributed couplings), unbounded control operators (case of boundary control), and in some situations if the diffusion operators are not the same.
The results of~\cite{alacras01, sicon03} have been recently extended by the author and L\'eautaud in~\cite{alaleaucont, alaleau11} to the case of partially coercive coupling operators (case of  localized couplings). The coupling coefficient  is assumed to be a sufficiently smooth function. The results are valid for localized as well as boundary control, and for the same diffusion operators. The geometric assumptions are that the control and coupling regions should both satisfy the Geometric Control Condition of Bardos Lebeau Rauch~\cite{blr92}
(see also
\cite{burq, burqgerard} for weaker smoothness assumptions on $\Omega$ and the coefficients of the elliptic operator $A$). This allows many situations for which the control and coupling regions do not necessarily meet. 

We can now turn to some consequences for the corresponding parabolic or Schr\"odinger systems.
The transmutation method
~\cite{seidman, miller, phung, EZ}, allows to deduce null controllability results for heat or Schr\"odinger equations from exact controllability results for the wave equation. This is an indirect method.  Thus, using the transmutation method, we prove in~\cite{alaleaucont, alaleau11},  null controllability for symmetric $2$-coupled systems of parabolic and diffusive equations. These results are valid in a multi-dimensional setting and under the condition that both the coupling and control regions satisfy the Geometric Control Condition. In particular, this allows many geometric situations for which these regions do not meet. A drawback of this indirect method is that these null controllability results are then, only valid under the usual geometric control conditions for the wave equations, whereas for instance null controllability for the scalar heat equation holds without geometric assumptions. 
\vskip 2mm

Let us now give an overview of the literature on insensitizing controls for the scalar wave equation and on controllability of $2$-coupled hyperbolic cascade systems.
As far as we know, the first results for insensitizing controls for the scalar wave equation are due to D\'ager~\cite{Dager06}. 
He proved the insensitizing boundary controllability, and the $\varepsilon$-insensitizing locally distributed controllability for the one-dimensional wave equation. In both situations, the control regions need not to meet the coupling region. Tebou~\cite{tebou08} has considered the same questions in the multi-dimensional framework. He partially extended D\'ager's results to the multi-dimensional wave equation for controls which are locally distributed in a region $\omega$, and furthermore coupling regions $O$ which necessarily meet $\omega$. More precisely, he proved the $\varepsilon$-insensitizing locally distributed controllability for arbitrary open subsets $\omega$ and $O$  such that $\omega \cap O \neq \emptyset$. He also proved the insensitizing locally distributed controllability under a strong geometric assumption, namely  that  both the control and coupling regions contain the same neighbourhood of a part $\Gamma_1$ of the boundary, that satisfies the usual multiplier condition. A slightly different analysis has been performed by Tebou in~\cite{tebou2011}. In this paper, the functional to be insensitized involves the trace of the normal derivative of the unknown on a part $\Gamma_1$ of the boundary. The existence of a locally distributed insensitizing control is proved under the strong geometric condition that $\Gamma_1$ satisfies the usual multiplier geometric condition and that the localized control region contains a neighbourhood of $\Gamma_1$. 

In a recent work, Rosier and de Teresa~\cite{RdT11} considered a $2$-coupled system of cascade hyperbolic equations under a strong hypothesis, that is a periodicity assumption of the semigroup associated to a single uncoupled equation. They give applications to $2$-coupled systems of cascade one-dimensional heat equations and to $2$-coupled systems of cascade Schr\"odinger equations in a $n$-dimensional interval with empty intersection between the control and coupling regions. The coupling coefficient is a characteristic function and thus is not a smooth function. Their method is linked to D\'ager's~\cite{Dager06} approach, based on the periodicity assumption of the semigroup for the single free equation. 
 Dehman L\'eautaud Le Rousseau~\cite{these-leautaud} (see also~\cite{DLRL}) consider coupled cascade wave systems
in a $\mathcal{C}^{\infty}$ compact connected riemannian manifold without boundary in the case of  locally distributed observations. Using a contradiction argument, they prove an observability inequality  for such systems, 
and further give the characterization of the minimal control time. Their proof uses micro-local analysis and the principle to work in weakened energy spaces for the unobserved component (see~\cite{alacras01, sicon03}). They deduce the corresponding controllability result thanks to the HUM method.
 
We also would like to refer to some books on coupled systems and on
control theory. In particular, the interested reader can find a series of results on the control of coupled systems issued from mechanical applications, and in particular for acoustic models, in Lasiecka~\cite{lasieckabk02}.  We refer to J.-L. Lions~\cite{lions} for the HUM method, together with examples of partial control of coupled systems in the case of infinite dimensional control systems and to Coron~\cite{coronbk07} for a nice introduction of the HUM method in a finite dimensional setting, a general approach of nonlinear control and various examples of nonlinear control for PDE's.

\vskip 2mm

{\bf Main results of the paper}

\vskip 2mm

We give a necessary and sufficient condition for the observability of $2$-coupled abstract hyperbolic cascade systems by a single observation. The observatility operator can be a bounded as well as an unbounded operator. The coupling operator is assumed to be only partially coercive as in~\cite{alaleau11}.  In the case of domains such that their boundaries have no contact of infinite order with its tangent (or for analytic boundaries),
this necessary and sufficient condition states that the support of the observation and coupling functions should satisfy (GCC). In the one-dimensional case, this shows that it is possible to drive back to equilibrium the solution of an uncontrolled wave equation locally coupled to a controlled scalar wave equation, the coupling being active on any non-empty open subset of arbitrary measure. We prove these results, by adapting the two-level energy method~\cite{alacras01, sicon03, alaleau11} to cascade systems. This method has been originally introduced to handle symmetric coupled systems for which the total energy of the solutions is conserved through time, whereas this property is lost for cascade systems.  We deduce controllability results for coupled wave systems. Using transmutation, we prove  null-controllability results for coupled heat or Schr\"odinger systems under geometric conditions. These conditions are those of hyperbolic systems. We give several applications of these results to $2$-coupled cascade wave, heat and Schr\"odinger systems. They can also be applied to higher order systems such as Petrowsky equations, or wave equations with variable coefficients. The main point of this paper is that these results are valid in a multi-dimensional framework, for locally distributed as well as boundary controls (resp. observations), and for localized couplings in situations for which the control/observations regions do not meet the localized coupling regions.
 They answer positively (and partially due to geometric conditions) to an open question raised by Luz de Teresa~\cite{DeT00} for forward cascade heat systems. 
 
 We further give applications to the question of existence of exact insensitizing locally distributed and boundary controls for the wave equation as introduced by J.-L. Lions~\cite{lions89}, obtaining new and complete results on this question, generalizing those of D\'ager~\cite{Dager06} and Tebou~\cite{tebou08} (see also~\cite{tebou2011}). This result is important since it allows to build exact controls for the scalar wave equation which are robust, since these controls both drive back the solution to equilibrium and insensitize a weighted observation of the solution, making this observation insensitive to small unknown perturbations of the initial data.
 
 Parts of these results (without proofs) have been announced in~\cite{arxiv}. It should be noted that it is not possible for insensitizing control of the heat equation to use results from coupled cascade wave systems through the transmutation method. This comes from the fact the $2$-coupled cascade systems issued from insensitizing control for the heat equation, are two coupled heat equations one being forward in time, the second being backward in time.

\section{Controllability and observability of $2$-coupled cascade hyperbolic systems by a single control/observation}
\subsection{Observability of $2$-coupled cascade hyperbolic systems by a single observation}
\setcounter{equation}{0}
We consider the following coupled cascade
hyperbolic system of order $2$

\begin{equation}\label{NSH}
\begin{cases}
u_1^{\prime\prime} + A u_1 = 0 \,,\\
u_2^{\prime\prime} + A u_2+ C_{21}u_1  = 0 \,,\\
(u_i,u_i^{\prime})(0)=(u_i^0,u_i^1) \mbox{ for }
i=1,2 \,,
\end{cases}
\end{equation}
\noindent where $H$ is an Hilbert space with norm $|\cdot|$ and scalar product $\langle \,,\,\rangle$ and $C_{21}$
is a bounded operator in $H$. We assume that $A$ satisfies
\begin{equation}\label{A1}
(A1)\ 
\begin{cases}
 A : D(A) \subset H \mapsto H \,, A^{\ast}=A\,,\\
\exists \; \omega>0\,,  |A u| \ge \omega |u|  \quad \forall \ u \in D(A) \,,\\
A \mbox{  has a compact resolvent}\,.
\end{cases}
\end{equation}
\noindent One can note that this system is a lower triangular system, that is it involves
a lower triangular operator when written as a second order equation in vectorial form. More precisely, we have

\begin{equation*}
\begin{pmatrix}
u_1\\
u_2
\end{pmatrix} 
^{\prime\prime}=
\begin{pmatrix}
-A &  0\\
-C_{21} & -A
\end{pmatrix}
\begin{pmatrix}
u_1\\
u_2
\end{pmatrix} 
\,.
\end{equation*}
\noindent We set $H_k=D(A^{k/2})$ for $k \in \N$, with the convention
$H_0=H$. The set $H_k$ is equipped with the norm $|\cdot|_k$
defined by $|A^{k/2} \cdot|$ and the associated scalar product. It
is a Hilbert space. We denote by $H_{-k}$ the dual space of
$H_k$ with the pivot space $H$. We equip $H_{-k}$ with the norm
$|\cdot|_{-k}=|A^{-k/2} \cdot|$. We define the energy space associated
to \eqref{NSH} by $\mathcal{H}=H_1^2\times
H_0^2$. 

The system \eqref{NSH} can then be reformulated as the first order abstract system

\begin{equation}\label{ANSH}
\begin{cases}
U^{\prime} = \mathcal{A} U \,,\quad \\
U(0)=U^0=(u_1^0,u_2^0,u_1^1,u_2^1)\,,
\end{cases}
\end{equation}
\noindent where $U=(u_1,u_2,v_1,v_2)$ and
$\mathcal{A}$ is the unbounded operator
in $\mathcal{H}$ with domain
$D(\mathcal{A})=H_2^2\times H_1^2$ defined by

\begin{equation}\label{A}
\mathcal{A}U=(v_1,v_2,-Au_1,-Au_2 - C_{21}u_1) \,.
\end{equation}
\noindent  Using semigroup theory, it is easy to establish the well-posedness of the abstract system
\eqref{ANSH} for initial data $U_0 \in
\mathcal{H}$. Moreover for initial data $U^0 \in D(\mathcal{A})$, 
we easily deduce that $Cu_1 \in \mathcal{C}^1([0,T]; H)$. Using classical results
on inhomogeneous Cauchy problems for the first order equation satisfied by $U_2=(u_2,u_2^{\prime})$ and the regularity of $Cu_1$, one can check that  for $U^0 \in D(\mathcal{A})$, the solution $U$ of \eqref{ANSH} is such that
$U \in \mathcal{C}([0,T]; D(\mathcal{A}))\cap  \mathcal{C}^1([0,T]; \mathcal{H})$.
Moreover, assuming that
$C_{21} 
\in \mathcal{L}(H_{k-1})$ for $k \in \mathbb{Z}^{\ast}$, the problem \eqref{ANSH} (and similarly \eqref{NSH}) is well-posed in $H_k^2 \times H_{k-1}^2$, that is
 if the initial data are in $H_k^2 \times H_{k-1}^2$, then  the solution
$U$ of \eqref{ANSH} (and similarly that of \eqref{NSH}) is  in 
$\mathcal{C}([0,T]; H_k^2 \times H_{k-1}^2)$. For a solution
$U=(u_1,u_2,v_1,v_2)$ of \eqref{NSH}, we have
$v_i=u_i^{\prime}$ for $i=1,2$. 
In the sequel, we will need several levels of energy of solutions $U$
of \eqref{ANSH}. For this, it is convenient to introduce some
further notation. For a solution
$U=(u_1, u_2,u_1^{\prime}, u_2^{\prime})$ of \eqref{ANSH}, we set
\begin{equation}\label{Ui}
U_i=(u_i,u_i^{\prime}) \mbox{ for } i=1,2\,.
\end{equation}
\noindent For $U_i \in H_k \times H_{k-1}$, we define the {\it local} energies of level $k$ as
\begin{equation}\label{enkU}
e_k(U_i)(t)=\frac{1}{2} \Big(
|A^{k/2}u_i|^2 + |A^{(k-1)/2}u_i^{\prime}|^2\Big) \,, \ 
k \in \mathbb{Z} \,, i=1,2\,.
\end{equation}
\noindent For $k=1$, we recover the natural energy of each component
of the state. For $k <1$, these energies are weakened energies.
We will also use this notation for more
general set of vector-valued functions $t \mapsto V(t)=(v_1(t),v_2(t))
\in H_k\times H_{k-1}$ for convenience, that is we define $e_k(V)(t)$
as above without further recalling this in the sequel.

We will also need {\it global} energies of level $k$. For this, we need to invert $\mathcal{A}$ on the
set of solutions $U=(u_1,u_2,v_1,v_2)$ of \eqref{ANSH}. We check that this is possible in the following proposition.
\begin{Proposition}\label{invertA}
Assume that $A$ satisfies $(A1)$ and define $\mathcal{A}$
as in \eqref{A}. Then $\mathcal{A}$ is invertible from
$D(\mathcal{A})$ on $\mathcal{H}$. Moreover
for any solution $U=(u_1,u_2,v_1,v_2)$ of \eqref{ANSH}, the equation
\begin{equation}\label{AW}
\mathcal{A}W=U \,,
\end{equation}
\noindent admits a unique solution $W=(w_1,w_2,r_1,r_2)$ given by
\begin{equation}\label{W}
\begin{cases}
w_1=-A^{-1}u_1^{\prime} \,,\\
w_2=-A^{-1}u_2^{\prime} + A^{-1}C_{21}A^{-1}u_1^{\prime} \,,\\
r_1=w_1^{\prime}=u_1 \,, r_2=w_2^{\prime}=u_2 \,.
\end{cases}
\end{equation}
\noindent Also, $W$ is then the solution of \eqref{ANSH}, associated
to the initial data $W(0)=(A^{-1}u_1^1,\linebreak -A^{-1}u_2^1 + A^{-1}C_{21}A^{-1}u_1^1,u_1^0,u_2^0)$.
\end{Proposition}
\begin{Proof}
We associate to $U=(u_1,u_2,v_1,v_2) \in \mathcal{H}$, the
vector $W=(w_1,w_2,r_1,r_2)$ defined by

\begin{equation}\label{Wgen}
\begin{cases}
w_1=-A^{-1}v_1 \,,\\
w_2=-A^{-1}v_2 + A^{-1}C_{21}A^{-1}v_1 \,,\\
r_1=u_1 \,, r_2=u_2 \,.
\end{cases}
\end{equation}
\noindent One can check that $W \in D(\mathcal{A})$ and satisfies \eqref{AW}. Hence $\mathcal{A}$ is
invertible. On the other hand $U$ satisfies \eqref{ANSH}. Thus applying the operator $\mathcal{A}^{-1}$ on both sides of the first
equation of \eqref{ANSH}, we deduce that $W$ also satisfies 
$$
W^{\prime} = \mathcal{A} W \,.
$$
\noindent Therefore, by definition of $\mathcal{A}$, we have $w_i^{\prime}=r_i$ for $i=1,2$. In a similar way, we have
$u_i^{\prime}=v_i$ for $i=1,2$. This, together with
\eqref{Wgen} imply that  $W$ is the unique solution of \eqref{W}. Moreover if $U$
is a solution of \eqref{ANSH}, then $W$ is the solution of \eqref{ANSH}
corresponding to initial data $W(0)=(A^{-1}u_1^1,-A^{-1}u_2^1 + A^{-1}C_{21}A^{-1}u_1^1,u_1^0,u_2^0)$, and given by \eqref{W}.
\end{Proof}
We have the following corollary, which proof is left to the
reader.

\begin{Corollary}\label{coroinduc}
We assume the hypotheses of Proposition~\ref{invertA}.
Then the equation 
\begin{equation}\label{eqWk}
\mathcal{A}^kW^k=U
\end{equation}
\noindent admits a unique solution $W^k \in D(\mathcal{A}^k)$ defined
by induction as 
\begin{equation}\label{inducW}
W^0=U \,, W^{k+1}=\mathcal{A}^{-1}W^k \,, k \in \mathbb{N}\,.
\end{equation}
\end{Corollary}

\begin{Remark}
In the above corollary, we can only assert that $W^k$ is
in $D(\mathcal{A}^k)$. If in addition, the bounded operator
$C_{2 1}$ satisfies $C_{21}H_{k-1} \subset H_{k-1}$ for $k \in \mathbb{N}^{\ast}$,
then $W^k \in H_{k+1}^2 \times H_k^2$.
\end{Remark}

\subsubsection{Main results for the observability of $2$-coupled cascade systems}
Our purpose is to establish indirect observability estimates for the system \eqref{NSH}. 
For this, we shall assume the following hypotheses $(A2)-(A5)$.

\medskip

We assume that the coupling operator $C_{2 1}$  satisfies
\begin{equation}\label{hypC}
(A2)
\begin{cases}
\, C_{21}^{\ast} \in \mathcal{L}(H_k) \mbox{ for } k \in \{0,1\} \,,\\
 \, ||C_{21}||=\beta\,,
|C_{21}w|^2 \le \beta \langle C_{21}w\,,w\rangle \quad \forall \ w \in H \,,\\
\exists \alpha>0 \mbox{ such that, }
\alpha\, |\Pi w|^2 \le \langle C_{21}w,w\rangle \quad
\forall \ w \in H \,,
\end{cases}
\end{equation}
\noindent where the operator $\Pi$ satisfies the assumptions
\begin{equation}\label{observabilityPi}
(A3)
\begin{cases}
\Pi \in \mathcal{L}(H)\,, \exists \ T_0>0, \forall \ T>T_0 \,, \exists \ C_2(T)>0 \mbox{ such that }\\
\forall \ (w^0,w^1) \in H_1 \times H
\mbox{ the solution } w \mbox{ of }\\
w'' + A w = 0 \,, (w,w^{\prime})(0)=  (w^0,w^1)   
\mbox{ satisfies }  \\
\int_0^T |\Pi w^{\prime}|^2 dt \geq C_2(T)e_1(W)(0) \,,
\end{cases}
\end{equation}
\noindent where $W=(w,w^{\prime})$. We denote by $G$ a given Hilbert space with norm $||\ ||_{G}$ and
scalar product $\langle \,, \rangle_{G}$. The space
$G$ will be identified to its dual space in all the sequel.
We make the following assumptions on the observability operator $\mathcal{\mathbf{B}}^{\ast}$  (the dual operator of the control operator $\mathcal{\mathbf{B}}$). 

We shall first assume that $\mathcal{\mathbf{B}}^{\ast}$  is an admissible observation operator for one equation, that is
\begin{equation}\label{admissibility}
(A4)
\begin{cases}
\mathcal{\mathbf{B}}^{\ast} \in \mathcal{L}(H_2\times H;G),\\
\forall \ T > 0 \ \exists \ C >0 , \mbox{ such that for all } (w^0,w^1) \in H_1 \times H \mbox{ and all } f  \in L^2([0,T];H)\,,\\
\mbox{the solution } w \mbox{ of }
w'' + A w = f \,, (w,w^{\prime})(0)=  (w^0,w^1)
\mbox{ satisfies }  \\
\int_0^T \| \mathcal{\mathbf{B}}^{\ast}(w,w^{\prime}) \|_G^2 dt \leq C \left( e_1(W)(0) + e_1(W)(T) + \int_0^T e_1(W)(t) dt + \int_0^T |f|^2 dt \right),
\end{cases}
\end{equation}
\noindent where $W=(w,w^{\prime})$. 

We further assume the following observability inequality for a single equation
\begin{equation}\label{observabilityB}
(A5)
\begin{cases}
\exists \ T_0>0, \forall \ T>T_0 \,, \exists \ C_1(T)>0 \mbox{ such that }\\
\forall \  (w^0,w^1) \in H_1 \times H \,,
\mbox{ the solution } w \mbox{ of }\\
w'' + A w = 0   \,,     (w,w^{\prime})(0)=  (w^0,w^1)
\mbox{ satisfies } \\
\int_0^T \| \mathcal{\mathbf{B}}^{\ast}(w,w^{\prime}) \|_G^2 dt \geq C_1(T)e_1(W)(0)\,. 
\end{cases}
\end{equation}

\begin{Remark}
The minimal times for which the two observability inequalities hold in $(A3)$ and $(A5)$ are not necessarily the same for the two observability operators, but here, we only consider the sup of these two optimal times to avoid too many notation.
\end{Remark}
 \begin{Lemma}\label{admiss}(Admissibility property)
 Assume the hypotheses $(A1)$, $(A4)$ and that $C_{21} \in \mathcal{L}(H)$, then for all $T>0$, there exists a constant $C=C(T)>0$ such that for all initial data $U^0 \in \mathcal{H}$, the solution of \eqref{NSH} satisfies the following direct inequality
 
 \begin{equation}\label{admissineq}
 \int_0^T ||\mathcal{\mathbf{B}}^{\ast}U_2||_G^2 \,dt \le C \Big(e_0(U_1)(0) +e_1(U_2)(0)\Big) \,.
 \end{equation}
  \end{Lemma}
  \begin{Remark}
  This Lemma establishes a hidden regularity property of the solutions, namely that for all $U^0 \in
  \mathcal{H}$, $\mathcal{\mathbf{B}}^{\ast}U_2 \in L^2([0,T];G)$.
  \end{Remark}
 \begin{Theorem}[Sufficient conditions]\label{NSobs}
Assume the hypotheses $(A1)-(A5)$. 
Then there exists $T_3>0$ such that for all $T>T_3$,
and all initial data $U^0 \in \mathcal{H}$,
the solution of \eqref{NSH} satisfies
the observability estimates
\begin{equation}\label{obsNSH}
\begin{cases}
d_1(T) \int_0^T||\mathcal{\mathbf{B}}^{\ast}U_2||_G^2
\ge  e_0(U_1)(0)\,,\\
d_2(T)\int_0^T||\mathcal{\mathbf{B}}^{\ast}U_2||_G^2
\ge  e_1(U_2)(0) \,, 
\end{cases}
\end{equation}
\noindent where the constants $d_i(T)>0$
are obtained thanks to \eqref{obse10} (in Lemma~\ref{obse100})
and \eqref{obse200} (in Lemma~\ref{obse20}), depend on $T$ and satisfy for $T$ sufficiently large
\begin{equation}\label{d12}
d_1(T) \le \frac{K}{T^3} \,, \, d_2(T) \le \frac{K}{T} \,.
\end{equation}
\noindent Moreover, the following estimates also hold
\begin{equation}\label{inte1U2m}
\int_0^T e_1(U_2) \le k_2(T) \int_0^T||\mathcal{\mathbf{B}}^{\ast}U_2||_G^2 \,,
\end{equation}
\noindent and
\begin{equation}\label{C21}
\int_0^T \langle C_{2 1}u_1, u_1 \rangle  \le r_2(T) \int_0^T||\mathcal{\mathbf{B}}^{\ast}U_2||_G^2 \,,
\end{equation}
\noindent where for $T$ sufficiently large
\begin{equation}\label{k2}
k_2(T) \le K\,, r_2(T) \le K/T^2 \,.
\end{equation}
\noindent Here the notation $K$
 stands for positive constants which does not depend on,  $T$, but depends on $\alpha, \beta, \gamma_0$ in an explicit way. Moreover, one can choose $T_3$
in the following form $T_3=
\max(T_0,T_1,T_2)$, where $T_0$ is introduced in
$(A3)$, $(A5)$ and $T_1, T_2$ are defined in
\eqref{T1}. In addition, if $C_{21}$ satisfies $(A3)$, then one can choose $\Pi=C_{21}$ so that the above properties hold.
\end{Theorem}
We also prove that the above conditions are optimal in the following theorem.
\begin{Theorem}[Necessary conditions]\label{2NEC}
Assume the hypotheses $(A1)$ and $(A4)$, and that 
\begin{equation}\label{A2prime}
(A2)^{\prime}\begin{cases}
C_{21}^{\ast} \in \mathcal{L}(H_k) \mbox{ for } k \in \{0,1\} \,,\\
 \, ||C_{21}||=\beta\,,
|C_{21}w|^2 \le \beta \langle C_{21}w\,,w\rangle \quad \forall \ w \in H \,,\\
\end{cases}
\end{equation}
\noindent holds. Assume that $\Pi=C_{21}$ does not satisfy $(A3)$ or that $\mathcal{\mathbf{B}}$ does not satisfy $(A5)$. Then there does not exist $T_3>0$ such that for all $T>T_3$, 
the following property holds
\begin{equation}\label{OBSP}
(OBS) \begin{cases}
\exists \ C>0 \mbox{ such that } \forall \ U^0 \in \mathcal{H} \mbox{ the solution of \eqref{NSH} } \mbox{ satisfies} \quad \\
C (e_0(U_1)(0) + e_1(U_2)(0)) \le C \int_0^T||\mathcal{\mathbf{B}}^{\ast}U_2||_G^2\,dt \,.
\end{cases}
\end{equation}
\end{Theorem}
\begin{Corollary}\label{2CNS}
Assume $(A1)$ and $(A4)$ and that $\Pi=C_{21}$ satisfies $(A2)^{\prime}$. Then \eqref{OBSP} holds if and only if $(A3)$ and $(A5)$ hold.
\end{Corollary}

\subsubsection{Proofs of the main results for the observability of $2$-coupled systems}
\begin{Proof} of Lemma~\ref{admiss}.
 Thanks to assumption $(A4)$ applied to the second equation of  \eqref{NSH}, we have for all
 $T>0$ there exists $C(T)>0$ such that
  \begin{equation}\label{inter1}
\int_0^T || \mathcal{\mathbf{B}}^{\ast}U_2||_G^2 \,dt \le C(T)\Big(e_1(U_2)(0) + e_1(U_2)(T)+
\int_0^T e_1(U_2)(t)\,dt + \int_0^T |C_{21}u_1|^2\,dt
\Big)\,.
 \end{equation}
 \noindent On the other hand, the usual energy estimates yield
 $$
 e_1(U_2)(t) \le C\Big( e_1(U_2)(0) +T \int_0^T |C_{21}u_1|^2\,dt\Big)\,.
 $$
 \noindent Integrating this inequality between $0$ and $T$ and using
 it for $t=T$ in \eqref{inter1}, we obtain \eqref{admissineq}.
  \end{Proof}
The proof of Theorem~\ref{NSobs} follows the idea of the two-level energy method such as developed in \cite{sicon03}. It consists in using two level of energies, the natural one for the observed component of the state
and the weakened energy of the unobserved component of the state. Here, we no longer consider
symmetric conservative coupled hyperbolic systems such as considered in \cite{sicon03}, or more recently in \cite{alaleaucont}. For the two-level
energy method in the situation of \cite{sicon03}, a crucial property was assumed for the abstract system (and proved for the applicative examples), that is direct  and observability inequalities for a single equation with a source term, with constants which are {\em uniform} with respect to the length $T$ of the time interval $[0,T]$. This property was proved thanks to the multiplier method for the usual PDE's (wave, Petrowsky,\ldots) in \cite{sicon03} and extended  in \cite{alaleaucont} under an abstract form which allows the use of the optimal geometric conditions of Bardos Lebeau and Rauch~\cite{blr92}. We shall use this result in the sequel so we recall that it reads as follows. 
\begin{Lemma}[\cite{alaleau11}, [Lemma 3.3, pp. 14]\label{AL}
We assume the hypotheses $(A1)$, $(A4)$ and $(A5)$. Then, 
there exist constants $\eta_0>0$ and $\alpha_0>0$ such that for all $T>T_0$, and for any solution $P=(p,p^{\prime})$ of  the nonhomogeneous equation
\begin{equation}\label{eqNH}
p^{\prime\prime} +Ap=f  \in L^2([0,T];H)\,,
\end{equation}
\noindent  the following uniform observability estimate holds 
\begin{equation}\label{obsabs} 
\eta_0\int_0^T ||\mathcal{\mathbf{B}}^{\ast}P||^2 \,dt \ge 
\int_0^T e_1(P)(t)\,dt -
\alpha_0\int_0^T|f|^2 \,dt \,.
\end{equation}
\noindent In a similar way if $(A1)$ and $(A3)$ hold then there exist $\gamma_0>0$ and $\delta_0>0$ such that for all $T>T_0$, and for any solution $P=(p,p^{\prime})$ of  \eqref{eqNH}, the following uniform observability estimate holds 
\begin{equation}\label{obsabspi} 
\gamma_0\int_0^T |\Pi p^{\prime}|^2  \,dt \ge 
\int_0^T e_1(P)(t)\,dt -
\delta_0\int_0^T|f|^2 \,dt \,.
\end{equation}
\end{Lemma}
\begin{Remark}
Note that we give a slightly different presentation and use different notations for the constants than
in the original Lemma 3.3 in \cite{alaleau11} without loss of generality.
\end{Remark}

We deduce from this lemma the following corollary.
\begin{Corollary}\label{ALcoroll}
Assume that $(A1)$ and $(A3)$ hold, then there exists a constant $\gamma_0>0$ such that  for all $T>T_0$ and for all
solutions $V=(v,v^{\prime})$ of
\begin{equation}\label{eqH}
\begin{cases}
v^{\prime\prime} +Av=0 \,,\\
(v,v^{\prime})(0)=(v^0,v^1) \,,
\end{cases}
\end{equation}
\noindent the
following observability inequality holds
\begin{equation}\label{hypobsabs}
T e_1(V)(0) \le \gamma_0
\int_0^T |\Pi v^{\prime}|^2 \,dt \,.
\end{equation}
\end{Corollary}
\begin{Proof}
We use the inequality \eqref{obsabspi} with $f=0$ and the conservation of energy. This gives the desired result.
\end{Proof}
\begin{Remark}
It should be noted that in $(A3)$, the dependence of the constant $C_2(T)$ on $T$ is not known. Here we can
prove that this constant can be chosen under the form $C_2(T)=T/\gamma_0$. This is an important property for the two-level
energy method.
\end{Remark}

{\bf Proof of Theorem~\ref{NSobs}:}

The proof will be divided in several Lemma. 

In the sequel, we will assume that the initial data for \eqref{ANSH} are in $D(\mathcal{A})$.
The final result for all initial data follows easily from the density of $D(\mathcal{A})$
in $\mathcal{H}$. We will also use the above notation without further specifying it. Furthermore, we set $C=C_{21}$ in the sequel of this subsection.

As in \cite{sicon03}, we obtain a first estimate using the coupling operator and involving the natural energy of $u_2$ and the weakened energy of $u_1$.

\begin{Lemma}\label{estcoup}
We assume $(A1)$. Let $U^0 \in D(\mathcal{A})$. Then
the solution of
\eqref{NSH} satisfies the estimate
\begin{equation}\label{eq10}
\int_0^T \langle C u_1,u_1\rangle \le
2 \eta e_0(U_1)(0) + \frac{1}{\eta}(
e_1(U_2)(T)+e_1(U_2)(0)) \,, \forall \ \eta>0 \,.
\end{equation}
\end{Lemma}

\begin{Proof} 
Since $(u_1,u_2)$ is a solution of \eqref{NSH},
we have

$$
\int_0^T \langle u_1^{\prime\prime}+Au_1,u_2\rangle
- \langle u_2^{\prime\prime}+Au_2 +C u_1 ,u_1\rangle=0\,,
$$
\noindent so that
\begin{equation}\label{est1}
\int_0^T \langle Cu_1,u_1 \rangle =
\Big[\langle u_1^{\prime},u_2\rangle
- \langle u_2^{\prime},u_1\rangle \Big]_0^T\,.
\end{equation}
\noindent We estimate the right hand side of
\eqref{est1} as follows. We set
$$
I(t)=\langle u_1^{\prime},u_2\rangle
- \langle u_2^{\prime},u_1\rangle \,.
$$
\noindent Then we have
$$
|I(t)| \le \frac{\eta}{2}\Big(|A^{-1/2}u_1^{\prime}|^2 +|u_1|^2\Big) +
\frac{1}{2\eta}\Big(|A^{1/2}u_2|^2 +|u_2^{\prime}|^2\Big) \,.
$$
\noindent Since the weakened energy $e_0(U_1)$ is conserved, we obtain the desired estimate.
\end{Proof}

\begin{Corollary}\label{e10}
We assume $(A1)-(A5)$. Let $U^0 \in D(\mathcal{A})$. Then
the solution of
\eqref{NSH} satisfies

\begin{equation}\label{eq13}
\int_0^T \langle Cu_1,u_1 \rangle \le
\frac{8 \gamma_0}{\alpha \, T} (e_1(U_2)(T)+e_1(U_2)(0)) \,.
\end{equation}
\end{Corollary}

\begin{Proof}
We define $W=(w_1,w_2,r_1,r_2)$ by \eqref{AW}. Then we proved in Proposition~\ref{invertA} that $W$ solves \eqref{W} and
\eqref{ANSH}, so that in particular $(w_1,w_1^{\prime})$ satisfies \eqref{eqH}. Thanks to the hypothesis $(A3)$ for $(w_1,w_1^{\prime})$
and to Lemma~\ref{AL}, \eqref{hypobsabs} holds for $V=W_1$.
We note that $w_1^{\prime}=u_1$. Thus,
thanks to the hypotheses $(A2)-(A3)$ and to \eqref{hypobsabs} together with \eqref{eq10} with the choice
$\eta=\frac{T\alpha}{4 \gamma_0}$, we obtain

$$ 
e_0(U_1)(0) \le \frac{8 \gamma_0^2}{\alpha^2 T^2} \big(e_1(U_2)(T)+ e_1(U_2)(0)\big) \,.
$$
\noindent Using this estimate in \eqref{eq10} with the above choice of $\eta$, we get \eqref{eq13}.
\end{Proof}

\begin{Lemma}\label{e20T}
Assume the hypotheses of Corollary~\eqref{e10}. Then,

\begin{equation}\label{eq15}
(e_1(U_2)(T)+e_1(U_2)(0)) \le
c_1 e_1(U_2)(0) +
\frac{c_2 \beta \gamma_0}{\alpha\, T} 
\int_0^Te_ 1(U_2)\,,
\end{equation}
\noindent and,

\begin{equation}\label{eq16}
\int_0^T \langle Cu_1,u_1 \rangle \le
\frac{c_3\gamma_0}{\alpha \, T} e_1(U_2)(0) +
\frac{c_4 \beta \gamma_0^2}{\alpha^2 \, T^2} 
\int_0^Te_ 1(U_2)\,.
\end{equation}
\end{Lemma}
\begin{Proof}
Let $(u_1,u_2)$ be a solution of \eqref{NSH}, then we have

$$
\langle u_2^{\prime\prime}+Au_2 +Cu_1,u_2^{\prime}\rangle =0 \,,
$$
\noindent so that

\begin{equation}\label{e2prime}
e_1^{\prime}(U_2)(t)=-\langle Cu_1,u_2^{\prime}\rangle(t) \,.
\end{equation}
\noindent Integrating this relation between
$0$ and $T$, we obtain
$$
e_1(U_2)(T)+ e_1(U_2)(0)\le 2 e_1(U_2)(0) + \frac{\alpha T}{16 \gamma_0} \int_0^T \langle Cu_1,u_1\rangle +
\frac{8\beta\gamma_0}{\alpha T}\int_0^T e_1(U_2) \,.
$$
\noindent We use \eqref{eq13} in this last estimate. This gives \eqref{eq15}. Using
\eqref{eq15} in \eqref{eq13}, we obtain \eqref{eq16}.
\end{Proof}

\begin{Lemma}\label{mininte2}
Assume the hypotheses of Corollary~\ref{e10}. Then,
\begin{equation}\label{eq18}
\int_0^T e_1(U_2)(t) \ge M\,T e_1(U_2)(0) \,,
\end{equation}
\noindent where $M$ is defined by \eqref{M}
and depends only on $\alpha,\beta,\gamma_0$.
\end{Lemma}

\begin{Proof}
We set 
\begin{equation}\label{a}
a=\frac{c_3 \beta \gamma_0}{2\alpha}\,,
\end{equation}
\noindent and
\begin{equation}\label{b}
b=\frac{c_4\beta^2 \gamma_0^2}{2\alpha^2}\,.
\end{equation}
\noindent We define
\begin{equation}\label{nu}
\nu=(a+ \sqrt{a^2+a+b})
\end{equation}
We integrate twice the two sides of \eqref{e2prime}, first between $0$ and $s$ and then between $0$ and $T$. This gives
$$
\int_0^T e_1(U_2)(t)=Te_1(U_2)(0)- \int_0^T (T-t) \langle Cu_1,u_2^{\prime}\rangle(t) \,.
$$
\noindent Using Young's inequality on the second term of the right hand side in the above relation, together with \eqref{hypC}
in assumption $(A2)$ and the definition of $e_1(U_2)$, we deduce that
$$
(1+\nu)\int_0^T e_1(U_2)(t) \ge Te_1(U_2)(0)
- \frac{\beta T^2}{2\nu} \int_0^T\langle Cu_1,u_1\rangle \,.
$$
\noindent Using \eqref{eq16} in this last estimate
and the definition of $a$, $b$ and $\nu$, we obtain \eqref{eq18} where
$M=M(\alpha,\beta,\gamma_0)$
is defined by

\begin{equation}\label{M}
M=\frac{\sqrt{a^2+a+b}}{(2a+1)\big(
a+\sqrt{a^2+a+b}\big)+a+2b} \,.
\end{equation}
\end{Proof}

\begin{Lemma}\label{obse20}
Assume the hypotheses of Theorem~\ref{NSobs}. We set
\begin{equation}\label{T1}
T_1=\frac{\sqrt{2c_4\alpha_0}\beta\gamma_0}{\alpha}
 \,,
T_2=\frac{\sqrt{2c_3\alpha_0\beta\gamma_0}}{
\sqrt{\alpha M}}\,,
T_3=\max\Big(T_0,T_1,T_2\Big)\,,
\end{equation}
\noindent where $T_0$ is introduced in $(A3)$ and $(A5)$. 
Then for all $T>T_3$, we have
\begin{equation}\label{obse200}
\eta_0 \int_0^T||\mathcal{\mathbf{B}}^{\ast}U_2||_G^2 \ge
\frac{M}{2T}\Big(T^2 - T_2^2\Big)e_1(U_2)(0)\,,
\end{equation}
\noindent and
\begin{equation}\label{inte1U2}
\int_0^Te_1(U_2) \le 2 \eta_0 \frac{T^2}{T^2 - T_2^2}
\int_0^T||\mathcal{\mathbf{B}}^{\ast}U_2||_G^2 \,.
\end{equation}

\end{Lemma}

\begin{Proof}
Applying Lemma~\ref{AL} for
the equation satisfied by $u_2$ in \eqref{NSH}, we deduce that \eqref{obsabs} holds, that is

\begin{equation}\label{obs2}
\eta_0 \int_0^T||\mathcal{\mathbf{B}}^{\ast}U_2||_G^2 \ge
\int_0^T e_1(U_2)(t) -\alpha_0\beta\int_0^T
\langle Cu_1,u_1\rangle \,.
\end{equation}
\noindent 

We use \eqref{eq16} and \eqref{eq18}
in this last inequality. Then for
all $T>T_3$, we obtain \eqref{obse200}. Using in a similar
way, \eqref{eq16} and \eqref{eq18} in \eqref{obs2}, together
with the definitions of $T_1$ and $T_2$, and
\eqref{obse200} in the resulting equation, we obtain \eqref{inte1U2}.
\end{Proof}

\noindent Hence, we proved that we can reconstruct
the initial data of the second component of
the state, which is coupled
to the first component, from the observation of this second component. We now have to prove that we can reconstruct the weakened energy of the first component from the observation of the second component.

\begin{Lemma}\label{obse100}
Assume the hypotheses of Theorem~\ref{NSobs}. We define $T_3$ as in \eqref{T1}. Then for all $T>T_3$, we have
\begin{equation}\label{obse10}
e_0(U_1)(0) \le
\frac{2 \eta_0 \gamma_0^2}{\alpha^2(T^2-T_2^2) T}\Big[\frac{c_3}{M} +
\frac{c_4\beta\gamma_0}
{\alpha}\Big]
\int_0^T||\mathcal{\mathbf{B}}^{\ast}U_2||_G^2 \,.
\end{equation}
\end{Lemma}

\begin{Proof}
Using \eqref{obs2} and \eqref{eq16}, and
since $T>T_3$, we obtain
\begin{equation}\label{Cu1}
\int_0^T\langle Cu_1,u_1\rangle \le
\Big(\frac{c_3\gamma_0}{\alpha T} \Big)e_1(U_2)(0) +
\frac{2c_4\eta_0\beta\gamma_0^2}{\alpha^2(T^2-T_2^2)}
\int_0^T||\mathcal{\mathbf{B}}^{\ast}U_2||_G^2 \,.
\end{equation}
\noindent We define $W=(w_1,w_2,r_1,r_2)$ by \eqref{AW}. Then we proved in Proposition~\ref{invertA} that $W$ solves \eqref{W} and
\eqref{ANSH}, so that in particular $(w_1,w_1^{\prime})$ satisfies \eqref{eqH}.  Hence thanks to the uniform observability inequality \eqref{hypobsabs}, and to $(A2)$, we deduce
that
\begin{equation}\label{e1tilde0}
e_0(U_1)(0) \le \frac{\gamma_0}{\alpha T}
\int_0^T\langle Cu_1,u_1\rangle \,.
\end{equation}

Inserting \eqref{obse200} in \eqref{Cu1}
and using \eqref{e1tilde0} we derive
\eqref{obse10}.
\end{Proof}

Thus, thanks to \eqref{obse200} and
\eqref{obse10},
we establish indirect observability for the triangular hyperbolic system \eqref{NSH}. We also deduce easily the estimates \eqref{inte1U2m} and \eqref{C21},
so that Theorem~\ref{NSobs} is proved.
\begin{Proof} of Theorem~\ref{2NEC}.
Assume first that $\mathcal{\mathbf{B}}^{\ast}$ does not satisfy $(A5)$. We argue by contradiction and assume that there exists $T_3>0$ such that $(OBS)$ holds for all $T>T_3$.
Then choosing the initial data $U^0$ under the form $=(0,u_2^0,0,u_2^1)$ where
$(u_2^0,u_2^1)$ is arbitrary in $H_1 \times H_0$ we have by uniqueness that $u_1=u_1^{\prime}
\equiv 0$
and $u_2$ is the solution of
\begin{equation}\label{U2O}
\begin{cases}
u_2^{\prime\prime} + A u_2=0 \,,\\
(u_2,u_2^{\prime})(0)=(u_2^0,u_2^1) \,,
\end{cases}
\end{equation}
\noindent whereas $(OBS)$ reduces to: there exists $T_3>0$ such that for all $T>T_3$ we have
$$
\begin{cases}
\exists \ C>0 \mbox{ such that } \forall \  (u_2^0,u_2^1)\in  H_1 \times H_0 \mbox{ the solution of \eqref{U2O} } \mbox{ satisfies} \quad \\
C (e_1(U_2)(0)) \le C \int_0^T||\mathcal{\mathbf{B}}^{\ast}U_2||_G^2\,dt \,.
\end{cases}
$$
\noindent Hence $\mathcal{\mathbf{B}}^{\ast}$ satisfies $(A5)$ which contradicts our hypothesis.

Assume now that $\Pi=C_{21}$ does not satisfy $(A3)$. We argue again by contradiction and assume that there exists $T_3>0$ such that $(OBS)$ holds for all $T>T_3$. We now choose
the initial data $U^0$ under the form $=(u_1^0,0,u_1^1,0)$ where
$(u_1^0,u_1^1)$ is arbitrary in $H_1 \times H_0$. Then $(OBS)$ reduces to: there exists $T_3>0$ such that for all $T>T_3$ we have
$$
\begin{cases}
\exists \ C_1>0 \mbox{ such that } \forall \  (u_1^0,u_1^1)\in  H_1 \times H_0 \mbox{ the solution of \eqref{NSH} } \mbox{ satisfies} \quad \\
C_1 e_0(U_1)(0) \le \int_0^T||\mathcal{\mathbf{B}}^{\ast}U_2||_G^2\,dt \,.
\end{cases}
$$
\noindent On the other hand using the admissibility assumption $(A4)$ together with
the usual energy estimates and the fact that $e_1(U_2)(0)=0$ we have
$$
\int_0^T||\mathcal{\mathbf{B}}^{\ast}U_2||_G^2\,dt \le C \int_0^T |C_{21}u_1|^2 \,dt \,.
$$
\noindent Hence there exists $C_2>0$ such that
$$
C_2 e_0(U_1)(0) \le \int_0^T |C_{21}u_1|^2 \,dt \quad \forall (u_1^0,u_1^1) \in H_1 \times H_0 \,.
$$
\noindent We set $w=-A^{-1}u_1^{\prime}$. Then we have $w^{\prime}=u_1$ and
$w^{\prime\prime}+ Aw=0$. We set $W=(w,w^{\prime})$. Then we have $(w,w^{\prime})(0)=(w^0,w^1)=(-A^{-1}u_1^1,u_1^0) \in H_2 \times
H_1$. 
Then, thanks to the above inequality, we have
for all $T>T_3$
\begin{equation}\label{DENSITY}
\int_0^T |C_{21}w^{\prime}|^2 \ge C_2 e_1(W)(0)  \, \quad \mbox{ for all the solutions of }
w^{\prime\prime} + Aw=0\,,
\end{equation}
\noindent with initial data $(w^0,w^1) \in H_2 \times H_1$. By density of $H_2 \times H_1$ in $H_1 \times H_0$, and continuity arguments,
we deduce that \eqref{DENSITY} holds for any $w$ solution of $w^{\prime\prime} + Aw=0$ with initial data
$(w^0,w^1) \in H_1 \times H_0$,  so that $C_{21}$ satisfies $(A3)$, which contradicts our hypothesis.
\end{Proof}
To handle the control problem, we shall
need to prove the admissibility and observability properties under a slightly different form
(mainly for the case $B \in \mathcal{L}(G,H)$). We have the following results.

 \begin{Lemma}\label{tech1}
 Assume the hypotheses of Proposition~\ref{invertA}.  Then there exist $C>0\,, C_1>0\,, C_2>0$ such that
 for all $W \in \mathcal{H}$,  the following properties hold for $Z=\mathcal{A}^{-1}W$ 
 
 \medskip
 
 \begin{description}
 \item (i) $e_0(Z_1)=e_{-1}(W_1)$,
 
 \item (ii) $e_0(W_2) \le C \Big( e_0(Z_1) + e_1(Z_2) \Big)$,
 
 \item (iii) $e_1(Z_2) \le C \Big(e_{-1}(W_1)+e_0(W_2)\Big)$,
 
 \item (iv) $C_1 \Big(e_0(Z_1) +e_1(Z_2)\Big) \le e_{-1}(W_1) + e_0(W_2) \le
 C_2 \Big( e_0(Z_1) + e_1(Z_2)\Big)$.
  \end{description}
 \end{Lemma}
 
 \begin{Proof}
 $Z=(z_1,z_2,z_1^{\prime} , z_2^{\prime})$ is defined as $z_1=-A^{-1}w_1^{\prime}$, $z_2=
 -A^{-1}w_2^{\prime} +A^{-1}C_{21}A^{-1}w_1^{\prime}$, $z_i^{\prime}=w_i$ for $i=1,2$.
 Therefore $(i)$ holds.  On the other hand, we have
 
 $$
 e_0(W_2)=\frac{1}{2}\Big(|w_2|^2+ |A^{-1/2}w_2^{\prime}|^2\Big) \le
 C\Big(e_0(Z_1)+ e_1(Z_2)\Big)\,.
$$
\noindent We also have

$$
e_1(Z_2)=\frac{1}{2}\Big(|A^{1/2}z_2|^2+ |z_2^{\prime}|^2\Big) \le
 C\Big(e_{-1}(W_1)+ e_0(W_2)\Big)\,.
$$
\noindent We deduce easily $(iv)$ thanks to $(i)-(ii)-(iii)$.
 \end{Proof}
 \begin{Remark}
 The operator $\mathcal{A}$ defined in \eqref{A} generates a $\mathcal{C}^0$-semigroup on $H_{-1}^2 \times H_{-2}^2$.  Hence due to the property of reversibility of time,
the Cauchy problem  $U^{\prime}=\mathcal{A} U$, $U(T)=U^T \in H_{-1}^2\times H_{-2}^2$  is well-posed, that is
has a unique solution in $\mathcal{C}^0([0,T];H_{-1}^2 \times H_{-2}^2)$. For the duality with the control problem, we shall need to work in different functional spaces which depend on the assumption on the control operator (case of bounded or unbounded control operator).

We set 
 \begin{equation}\label{X-1}
 X_{-1}=H_{-1}\times H \times H_{-2} \times H_{-1} \,.
 \end{equation}
 \noindent Since $X_{-1}
 \subset H_{-1}^2\times H_{-2}^2$, we can also solve the Cauchy problem
 with $U^T \in X_{-1}$.  In a similar way, we set
 \begin{equation}\label{X1}
 X_{1}=H \times H_1 \times H_{-1}\times H \,.
 \end{equation}
 \noindent Since $X_{1}
 \subset H_{-1}^2\times H_{-2}^2$, we can also solve the Cauchy problem
 with $U^T \in X_{1}$. 
  \end{Remark}

\begin{Lemma}\label{obsdir}
Assume $(A1)-(A5)$. Let $T>0$ be given.
For $W^T=(w_1^T, w_2^T, q_1^T, q_2^T) \in X_{-1}$, we denote by $W=(w_1,w_2, w_1^{\prime}, w_2^{\prime})$ the unique solution
in $\mathcal{C}^0([0,T]; H_{-1}^2 \times H_{-2}^2)$ of 
\begin{equation}\label{WT}
\begin{cases}
w_1^{\prime\prime} + Aw_1=0 \,, \\
w_2^{\prime\prime} + A w_2 + C_{21} w_1=0 \,, \\
W_{|t=T}=W^T \,.
\end{cases}
\end{equation}
\noindent Then $W$ satisfies the following properties
\begin{itemize}
\item (i) $W \in \mathcal{C}^0([0,T]; X_{-1})$,

\item (ii) There exists $C_1=C_1(T)>0$, such that

\begin{equation}\label{directweak1}
C_1 \int_0^T || \mathcal{\mathbf{B}}^{\ast}Z_2||_G^2 \,dt \le e_{-1}(W_1)(0) + e_0(W_2)(0)\,,
\end{equation}
\noindent where $Z=\mathcal{A}^{-1}W$.

\item (iii) For all $T>T_3$, where $T_3$ is given in Theorem~\ref{NSobs}, there exists $C_2=C_2(T)>0$ such that

\begin{equation}\label{directweak2}
e_{-1}(W_1)(0) \le C_2  \int_0^T || \mathcal{\mathbf{B}}^{\ast}Z_2||_G^2 \,dt \,,
\end{equation}
\begin{equation}\label{directweak3}
e_0(W_2)(0) \le C_2 \int_0^T || \mathcal{\mathbf{B}}^{\ast}Z_2||_G^2 \,dt \,.
\end{equation}

\item (iv) Assume furthermore that $\mathcal{\mathbf{B}}^{\ast}(w,w^{\prime})=B^{\ast}w^{\prime}$ where $B \in \mathcal{L}(G,H)$ is such that
$(A4)$ and $(A5)$ hold. Then
properties $(ii)-(iii)$ become

\begin{equation}\label{directweak1bis}
C_1 \int_0^T ||B^{\ast}w_2||_G^2 \,dt \le e_{-1}(W_1)(0) + e_0(W_2)(0)\,,
\end{equation}
\noindent and 

for all $T>T_3$
\begin{equation}\label{directweak2bis}
e_{-1}(W_1)(0) \le C_2  \int_0^T || B^{\ast}w_2||_G^2 \,dt \,,
\end{equation}
\begin{equation}\label{directweak3bis}
e_0(W_2)(0) \le C_2 \int_0^T || B^{\ast}w_2||_G^2 \,dt \,,
\end{equation}
\noindent with the same constants $C_1$ and $C_2$ than in $(ii)-(iii)$.
\end{itemize}
\end{Lemma}
\begin{Remark}
This Lemma shows that $Z$ as above defined satisfies a hidden regularity result, namely that for all initial data $W^0 \in X_{-1}$, $\mathcal{\mathbf{B}}^{\ast}Z_2
\in L^2([0,T];G)$.
\end{Remark}
\begin{Proof}
Since $W \in \mathcal{C}^0([0,T];H_{-1}^2 \times H_{-2}^2)$, $w_1 \in
\mathcal{C}([0,T]; H_{-1})$. Thanks to assumption $(A2)$, $C_{21}^{\ast}  \in \mathcal{L}(H_{1})$, thus we have $C_{21} \in \mathcal{L}(H_{-1})$, thus
$w_2$ is a solution of

$$
\begin{cases}
w_2^{\prime\prime} +Aw_2=-C_{21}w_1 \in \mathcal{C}([0,T]; H_{-1})\,, \\
(w_{2})_{|t=T}=w_2^T \in H \,,\\
(w_{2}^{\prime})_{|t=T}=q_2^T \in H_{-1} \,,
\end{cases}
$$
\noindent so that $(w_2,w_2^{\prime}) \in \mathcal{C}([0,T]; H \times H_{-1})$ by uniqueness. This yields $W \in \mathcal{C}^0([0,T]; X_{-1})$.

We set $Z=\mathcal{A}^{-1}W$. Thanks to \eqref{admissineq} together with properties $(i)$ and $(iii)$ in Lemma~\ref{tech1}, we easily deduce \eqref{directweak1}. This proves $(ii)$.

Thanks to $(ii)$ in Lemma~\ref{tech1}, we have

$$
e_0(W_2)(0) \le C\Big(e_0(Z_1)(0)+ e_1(Z_2)(0)\Big) \,.
$$
\noindent This together with \eqref{obsNSH} yield

$$
e_0(W_2)(0) \le C\Big(d_1(T) + d_2(T)\Big) \int_0^T ||\mathcal{\mathbf{B}}^{\ast}Z_2||_G^2 \,dt  \,.
$$
\noindent On the other hand, thanks to $(i)$ in Lemma~\ref{tech1} and
to \eqref{obsNSH}, we have

$$
e_{-1}(W_1)(0) \le d_1(T) \int_0^T ||\mathcal{\mathbf{B}}^{\ast}Z_2||_G^2 \,dt  \,.
$$
\noindent Thus \eqref{directweak2}-\eqref{directweak3} hold with $C_2=\max\Big(C(d_1(T)+d_2(T)), d_1(T)\Big)$. This proves $(iii)$.

The properties $(iv)$ follow easily from the hypothesis on 
$\mathcal{\mathbf{B}}^{\ast}(w,w^{\prime})$ and from the definition of $Z$ which implies that
$z_2^{\prime}=w_2$.
\end{Proof}

\begin{Lemma}\label{obsdirunbounded}
Assume $(A1)-(A5)$. Let $T>0$ be given.
For $W^T=(w_1^T, w_2^T, q_1^T, q_2^T) \in X_{1}$, we denote by $W=(w_1,w_2, w_1^{\prime}, w_2^{\prime})$ the unique solution
in $\mathcal{C}^0([0,T]; H_{-1}^2\times H_{-2}^2)$ of  \eqref{WT}
Then $W \in \mathcal{C}^0([0,T]; X_{1})$.
\end{Lemma}
\begin{Proof}
The proof is similar to the proof of Lemma~\ref{obsdir} and is left to the reader.
\end{Proof}

\subsection{Controllability of $2$-coupled cascade hyperbolic systems by a single control}

We apply the HUM method~\cite{lions, komornik1997} to deduce from the indirect observability inequality obtained in the previous section, an indirect exact controllability result for the control problem. Prior to this, we will recall for the sake of completeness, the transposition method (see~\cite{lions, komornik1997}) which allows to define the solutions of the control problem.

\subsubsection{The transposition method for the scalar wave equation}
We consider the control problem
\begin{equation}\label{SGwave}
\begin{cases}
y^{\prime\prime} + A y= Bv\,,\\
(y,y^{\prime})(0)=(y^0,y^1) \,,
\end{cases}
\end{equation}
\noindent where $A$ satisfies $(A1)$. We consider two cases: 
\begin{itemize}
\item (i) either $B \in  \mathcal{L}(G;H)$ (bounded control operator). In this case, we define
$\mathcal{\mathbf{B}}^{\ast}(w,w^{\prime})=B^{\ast}w^{\prime}$. We also set $\mathcal{F}_{-1}=H \times H_{-1}$
and $\mathcal{F}^{\ast}_{-1}=H_1\times H$.

\item or (ii)  $B \in \mathcal{L}(G, H_2^{\prime})$ (unbounded control operator).  In this case, we define
$\mathcal{\mathbf{B}}^{\ast}(w,w^{\prime})=B^{\ast}w$. We also set $\mathcal{F}_1=H_1 \times H$
and $\mathcal{F}^{\ast}_1=H\times H_{-1}$.

 \end{itemize}
 We assume that $\mathcal{\mathbf{B}}^{\ast}$ satisfies the assumption $(A4)$. 
 
 \begin{Definition} 
 \begin{itemize} 
 
 \item Case $(i)$:
 Let us first assume that $\mathcal{\mathbf{B}}^{\ast}(w,w^{\prime})=B^{\ast}w^{\prime}$ with  $B \in  \mathcal{L}(G;H)$.
 Let $(y^0,y^1) \in \mathcal{F}^{\ast}_{-1}$ and $v \in L^2_{loc}([0,\infty);G)$ be fixed arbitrarily.
 We say that $(y,y^{\prime})$ is a solution by transposition of \eqref{SGwave} if $(y,y^{\prime}) \in
 \mathcal{C}([0,\infty); \mathcal{F}_{-1}^{\ast})$ satisfies
 \begin{multline}\label{transposed}
\langle y^{\prime}(T)\,, w^0_T\rangle_{H,H} -\langle y(T)\,, w^1_T\rangle_{H_1\,,H_{-1}}=
\langle y^1\,, w(0)\rangle_{H,H} -\langle y^0\,, w^{\prime}(0)\rangle_{H_1\,,H_{-1}} +\\
\int_0^T \langle v \,,
B^{\ast} w\rangle_{G} \,dt \quad \forall \ T>0\,, \forall \  (w^0_T, w^1_T) \in \mathcal{F}_{-1} \,,
 \end{multline}
 \noindent where $w$ is the solution of
 \begin{equation}\label{wavehomogene}
 \begin{cases}
w^{\prime\prime} + A w= 0\,,\\
(w,w^{\prime})(T)=(w^0_T, w^1_T) \,.
\end{cases}
 \end{equation}
 \item Case $(ii)$:
 Let us assume now that $\mathcal{\mathbf{B}}^{\ast}(w,w^{\prime})=B^{\ast}w$ where $B \in \mathcal{L}(G, H_2^{\prime})$.
  Let $(y^0,y^1) \in \mathcal{F}^{\ast}_{1}$ and $v \in L^2_{loc}([0,\infty);G)$ be fixed arbitrarily.
 We say that $(y,y^{\prime})$ is a solution by transposition of \eqref{SGwave} if $(y,y^{\prime}) \in
 \mathcal{C}([0,\infty); \mathcal{F}_{1}^{\ast})$ satisfies
 \begin{multline}\label{transposedu}
\langle y^{\prime}(T)\,, w^0_T\rangle_{H_{-1},H_1} -\langle y(T)\,, w^1_T\rangle_{H\,,H}=
\langle y^1\,, w(0)\rangle_{H_{-1},H_1} -\langle y^0\,, w^{\prime}(0)\rangle_{H\,,H} + \\
\int_0^T \langle v \,,
B^{\ast} w\rangle_{G} \,dt \quad \forall \ T >0\,, \forall \ (w^0_T, w^1_T) \in \mathcal{F}_{1} \,,
 \end{multline}
 \noindent where $w$ is the solution of \eqref{wavehomogene}.
\end{itemize}
 \end{Definition}
 \begin{Remark}
 It is well-known that thanks to $(A4)$, there exists a unique solution by transposition to
\eqref{SGwave} which depends continuously on the data 
$(y^0,y^1)$ and on $v$. Let us state give a sketch of the proof for the sake of completeness. The solution $w$ of \eqref{wavehomogene} is such that $(w(0), w^{\prime}(0))$ depends continuously
on $(w^0_T, w^1_T)\in \mathcal{F}_{-1}$ (resp.  on $(w^0_T, w^1_T)\in \mathcal{F}_{1}$) in the case $(i)$ (resp. in the case $(ii)$) . Let us assume that we are in the case $(i)$. Thanks to $(A4)$ with $f=0$ and applied to $Z=(z, z^{\prime})=(-A^{-1}w^{\prime}, w)$, we have
\begin{multline*}
\int_0^T\| \mathcal{\mathbf{B}}^{\ast}(z,z^{\prime}) \|_G^2 dt=
\int_0^T\| B^{\ast} w \|_G^2 dt 
 \leq \\
 C_T e_1(Z)(0) = C_T\| (w(0),w^{\prime}(0))\|_{H \times H_{-1}}^2 \le 
 C_T \| (w^0_T,w^{1}_T)\|_{H \times H_{-1}}^2 \,,
\end{multline*}
\noindent where $C_T$ is a generic constant which depends on $T$.
Hence the right hand side of \eqref{transposed} defines a continuous linear form with respect to
$(w^0_T, w^1_T) \in \mathcal{F}_{-1}$, and moreover this linear form depends continuously on $T$ for all $T>0$.
This implies that for all $T>0$ there exists a unique  solution to \eqref{transposed} and that this solution depends continuously
on $T$. Let us now assume that we are in the case $(ii)$. Then, thanks to $(A4)$ with $f=0$ and applied to $W$, we have
\begin{multline*}
\int_0^T\| \mathcal{\mathbf{B}}^{\ast}(w,w^{\prime}) \|_G^2 dt=
\int_0^T\| B^{\ast} w \|_G^2 dt 
 \leq \\
 C_T e_1(W)(0) = C_T\| (w(0),w^{\prime}(0))\|_{H_1 \times H}^2 \le 
 C_T \| (w^0_T,w^{1}_T)\|_{H_1 \times H}^2\,,
\end{multline*}
\noindent where $C_T$ is a generic constant which depends on $T$. Thus, the right hand side of \eqref{transposedu} defines a continuous linear form with respect to
$(w^0_T, w^1_T) \in \mathcal{F}_{1}$, and moreover this linear form depends continuously on $T$ for all $T>0$.
This implies that for all $T>0$ there exists a unique  solution to \eqref{transposedu} and that this solution depends continuously
on $T$.
\end{Remark}
\begin{Remark}
We can, without loss of generality, reverse the time, changing $t$ in $T-t$ in the dual problem \eqref{wavehomogene} in the above definition, that is we can fix the initial
data instead of fixing the final data (with the appropriate minor changes).
\end{Remark}
 \begin{Remark}\label{INH}
Note also that we can equivalently replace $0$ on the right hand side of \eqref{wavehomogene} by any $f$
in $L^2((0,T);H_{-1})$ in the case $(i)$ (resp. in $L^2((0,T);H)$ in the case $(ii)$) for any $T>0$. More precisely,
in the case $(i)$, we can define a solution by transposition by setting
\begin{multline}\label{transposedINH}
\langle y^{\prime}(T)\,, w^0_T\rangle_{H,H} -\langle y(T)\,, w^1_T\rangle_{H_1\,,H_{-1}}+\int_0^T \langle y\,,f \rangle_{H_1,H_{-1}}\,dt =
\langle y^1\,, w(0)\rangle_{H,H} -\\ \langle y^0\,, w^{\prime}(0)\rangle_{H_1\,,H_{-1}} +
\int_0^T \langle v \,,
B^{\ast} w\rangle_{G} \,dt \quad \forall \ T>0\,, \forall \  (w^0_T, w^1_T) \in \mathcal{F}_{-1}\,, \forall \ f
\in L^2((0,T);H_{-1})\,,
 \end{multline}
 \noindent where $w$ is the solution of
 \begin{equation}\label{wavehomogeneINH}
 \begin{cases}
w^{\prime\prime} + A w= f\,,\\
(w,w^{\prime})(T)=(w^0_T, w^1_T) \,.
\end{cases}
 \end{equation}
 \noindent In the case $(ii)$, we can define a solution by transposition by setting
\begin{multline}\label{transposeduINH}
\langle y^{\prime}(T)\,, w^0_T\rangle_{H_{-1},H_1} -\langle y(T)\,, w^1_T\rangle_{H\,,H}+\int_0^T \langle y\,,f \rangle_{H,H}\,dt =
\langle y^1\,, w(0)\rangle_{H_{-1},H_1} -\\ \langle y^0\,, w^{\prime}(0)\rangle_{H\,,H} +
\int_0^T \langle v \,,
B^{\ast} w\rangle_{G} \,dt \quad \forall \ T>0\,, \forall \  (w^0_T, w^1_T) \in \mathcal{F}_{1} \,, \forall \ f
\in L^2((0,T);H)\,,
 \end{multline}
 \noindent where $w$ is the solution of \eqref{wavehomogeneINH}.
 In both cases, one has to consider respectively the right hand sides of \eqref{transposedINH} and \eqref{transposeduINH} as continuous
 linear forms in the appropriate spaces, with respect to $(w^0_T,w^1_T,f)$.
 \end{Remark}

\subsubsection{The transposition method for cascade control systems}
We consider the control problem

\begin{equation}\label{CTH}
\begin{cases}
y_1^{\prime\prime} + A y_1 +C_{21}^{\ast}y_2= 0 \,,\\
y_2^{\prime\prime} + A y_2= Bv \,,\\
(y_i,y_i^{\prime})(0)=(y_i^0,y_i^1) \mbox{ for }
i=1,2 \,,
\end{cases}
\end{equation}
\noindent where either $B \in  \mathcal{L}(G;H)$ (bounded control operator) or
 $B \in \mathcal{L}(G, H_2^{\prime})$ (unbounded control operator). We can easily adapt the definition of solutions by transposition to this coupled 
 system in the sense below. 
 \begin{Definition}\label{transcascade}
 \begin{itemize}
\item (i) Let $\mathcal{\mathbf{B}}^{\ast}(w,w^{\prime})=B^{\ast}w^{\prime}$ 
with $B \in \mathcal{L}(G,H)$. We set 
\begin{equation}\label{Xstar-1}
X_{-1}^{\ast}=H_2\times H_1 \times H_1 \times H\,.
\end{equation}
\noindent We assume that $\mathcal{\mathbf{B}}^{\ast}$ satisfies $(A4)$. We say that $(y_1,y_2,y_1^{\prime}, y_2^{\prime}) \in \mathcal{C}([0,\infty);X_{-1}^{\ast})$ is a solution by transposition of \eqref{CTH} if 
\begin{multline}\label{TRANSi}
\int_0^T \langle v, B^{\ast}w_2\rangle_G \,dt=
\langle y_1^{\prime}(T), w_1(T)\rangle_{H_1,H_{-1}} -
\langle y_1(T), w_1^{\prime}(T)\rangle_{H_2,H_{-2}} +\\
\langle y_2^{\prime}(T), w_2(T)\rangle_{H,H} - 
\langle y_2(T), w_2^{\prime}(T)\rangle_{H_1,H_{-1}} - \\
\Big(
\langle y_1^1, w_1(0)\rangle_{H_1,H_{-1}} -
\langle y_1^0, w_1^{\prime}(0)\rangle_{H_2,H_{-2}} +\\
\langle y_2^1, w_2(0)\rangle_{H,H} - \langle y_2^0, w_2^{\prime}(0)\rangle_{H_1,H_{-1}}
\Big)
 \,, \forall \ T>0\,, \forall \ 
W^T \in X_{-1}\,,
\end{multline}
\noindent where $W=(w_1,w_2,w_1^{\prime},w_2^{\prime})$ is the solution of
\begin{equation}\label{tildeW}
\begin{cases}
w_1^{\prime\prime} + Aw_1=0 \,, \\
w_2^{\prime\prime} + A w_2 + C_{21} w_1=0 \,, \\
W_{|t=T}=W^T \,,
\end{cases}
\end{equation}
\noindent 
\item (ii) Let $\mathcal{\mathbf{B}}^{\ast}(w,w^{\prime})=B^{\ast}w$ with $B \in\mathcal{L}(G, H_2^{\prime})$. 
We set
\begin{equation}\label{Xstar1}
X_{1}^{\ast}=H_1\times H\times H \times H_{-1}\,.
\end{equation}
We assume that $\mathcal{\mathbf{B}}^{\ast}$ satisfies $(A4)$. We say that $(y_1,y_2,y_1^{\prime}, y_2^{\prime}) \in \mathcal{C}([0,\infty);X_{1}^{\ast})$ is a solution by transposition of \eqref{CTH}
if
\begin{multline}\label{TRANSu}
\int_0^T \langle v, B^{\ast}w_2\rangle_G \,dt=
\langle y_1^{\prime}(T), w_1(T)\rangle_{H,H} -
\langle y_1(T), w_1^{\prime}(T)\rangle_{H_1,H_{-1}} +\\
\langle y_2^{\prime}(T), w_2(T)\rangle_{H_{-1},H_1} - 
\langle y_2(T), w_2^{\prime}(T)\rangle_{H,H} - \\
\Big(
\langle y_1^1, w_1(0)\rangle_{H,H} -
\langle y_1^0, w_1^{\prime}(0)\rangle_{H_1,H_{-1}} +\\
\langle y_2^1, w_2(0)\rangle_{H_{-1},H_1} - 
\langle y_2^0, w_2^{\prime}(0)\rangle_{H,H}
\Big)
 \,, \forall \ T>0\,, \forall \ 
W^T \in X_{1}\,,
\end{multline}
\noindent where $W=(w_1,w_2,w_1^{\prime},w_2^{\prime})$ is the solution of
\eqref{tildeW}.
\end{itemize}
\end{Definition}
\begin{Remark}
The dual problem \eqref{tildeW} is not conservative as for the case of the dual problem for a single wave equation. However the
solution of the non homogeneous equation in $w_2$ depends continuously on the final data for $w_2$ and on the source term $w_1$, which
itself depends continuously on the final data for $w_1$. This implies the existence of a unique solution in the sense of transposition of the
above control problem in cascade, which has the desired regularity and depends continuously on the initial data and on the control $v$.
\end{Remark}
\begin{Remark}\label{rktrans}
Assume that $(y_1,y_2,y_1^{\prime}, y_2^{\prime})$ is a solution by transposition of \eqref{CTH} (in the required space depending on case $(i)$ or $(ii)$). We can choose final data $W^T$ such that the final data for $w_1$ are vanishing in \eqref{tildeW}. Hence $w_1 \equiv 0$ and $w_2$ is
any solution of the homogeneous scalar wave equation (in the appropriate space). We deduce then easily
that $(y_2, y_2^{\prime})$ is a solution by transposition
 of the scalar equation 
 \begin{equation}\label{Y2t}
 \begin{cases}
y_2^{\prime\prime} + A y_2= Bv\,,\\
(y_2,y_2^{\prime})(0)=(y_2^0,y_2^1) \,.
\end{cases}
\end{equation}
\noindent Using the Remark~\ref{INH} with $f=-C_{21}w_1$  for the definition of the solution by
transposition $y_2$ for \eqref{Y2t}, together with the definition~\ref{transcascade} for the solutions by transposition for
the cascade control system, we deduce the following properties
\begin{itemize}
\item (i) If $\mathcal{\mathbf{B}}^{\ast}(w,w^{\prime})=B^{\ast}w^{\prime}$ 
with $B \in \mathcal{L}(G,H)$, then $(y_1, y_1^{\prime})$ satisfies
\begin{multline}\label{TRANSy1}
\langle y_1^{\prime}(T), w_1(T)\rangle_{H_1,H_{-1}} -
\langle y_1(T), w_1^{\prime}(T)\rangle_{H_2,H_{-2}} -\\
\Big(
\langle y_1^1, w_1(0)\rangle_{H_1,H_{-1}} -
\langle y_1^0, w_1^{\prime}(0)\rangle_{H_2,H_{-2}}\Big)=\\
-\int_0^T \langle y_2, C_{21} w_1 \rangle_{H_1,H_{-1}}
 \,, \forall \ T>0\,, \forall \ 
W_1^T \in H_{-1}\times H_{-2}\,,
\end{multline}
\noindent where $W_1=(w_1,w_1^{\prime})$ is the solution of
\begin{equation}\label{W1}
\begin{cases}
w_1^{\prime\prime} + Aw_1=0 \,, \\
(w_1,w_1^{\prime})_{|t=T}=W_1^T=(w_{1,T}^0,  w_{1,T}^1)\,.
\end{cases}
\end{equation}
\noindent 
\item (ii) If $\mathcal{\mathbf{B}}^{\ast}(w,w^{\prime})=B^{\ast}w$ with $B \in\mathcal{L}(G, H_2^{\prime})$, then
$(y_1, y_1^{\prime})$ satisfies 
\begin{multline}\label{TRANSuy1}
\langle y_1^{\prime}(T), w_1(T)\rangle_{H,H} -
\langle y_1(T), w_1^{\prime}(T)\rangle_{H_1,H_{-1}} -\\
\Big(
\langle y_1^1, w_1(0)\rangle_{H,H} -
\langle y_1^0, w_1^{\prime}(0)\rangle_{H_1,H_{-1}}
\Big)=\\
-\int_0^T \langle y_2, C_{21} w_1 \rangle_{H,H}
 \,, \forall \ T>0\,, \forall \ 
W_1^T \in H \times H_{-1}\,,
\end{multline}
\noindent where $W_1=(w_1,w_1^{\prime})$ is the solution of
\eqref{W1}.
\end{itemize}
\noindent The converse is also true. If $y_2$ solves \eqref{Y2t} and $y_1$ satisfies \eqref{TRANSy1} (resp.
\eqref{TRANSuy1}) in the case $(i)$ (resp. $(ii)$), then $(y_1,y_2, y_1^{\prime}, y_2^{\prime})$ is the corresponding solution by transposition of the cascade
system \eqref{CTH}.
\end{Remark}

\subsubsection{Controllability results for cascade systems}
\begin{Theorem}\label{control2}
Assume the hypotheses $(A1)-(A5)$. We define $T_3>0$ as in Theorem~\ref{NSobs}. We have the following properties.
\begin{itemize}
\item (i) Let $\mathcal{\mathbf{B}}^{\ast}(w,w^{\prime})=B^{\ast}w^{\prime}$ 
where $B \in \mathcal{L}(G,H)$ is such that $(A4)-(A5)$ holds. Then, for  all $T >T_3$ ,
and all $Y_0 \in X_{-1}^{\ast}$, there exists a control function $v \in L^2((0,T);G)$ such that the solution $Y=(y_1,y_2,y_1^{\prime},y_2^{\prime})$ of \eqref{CTH} satisfies $Y(T)=0$.
\item (ii) Let $\mathcal{\mathbf{B}}^{\ast}(w,w^{\prime})=B^{\ast}w$ where $B \in\mathcal{L}(G, H_2^{\prime})$ is such that $(A4)-(A5)$ holds.  
Then, for all $T >T_3$,
and all $Y_0 \in X_{1}^{\ast}$, there exists a control function $v \in L^2((0,T);G)$ such that the solution $Y=(y_1,y_2,y_1^{\prime},y_2^{\prime})$ of \eqref{CTH} satisfies $Y(T)=0$.
\end{itemize}
\end{Theorem}
\begin{Proof}
We first consider the case $(i)$.  Let $Y^0=(y_1^0,y_2^0,y_1^1,y_2^1) \in X_{-1}^{\ast}$. We consider the bilinear form $\Lambda$ on $X_{-1} $ defined by

\begin{equation}\label{Lambda}
\Lambda(W^T,\widetilde{W}^T)= \int_0^T \langle B^{\ast}w_2, B^{\ast}\widetilde{w_2}\rangle_G\,dt \,,
\forall \ W^T, \widetilde{W}^T \in X_{-1} \,,
\end{equation}
\noindent and the linear form on $X_{-1}$ defined for all $W^T  \in X_{-1}$
by

\begin{equation}\label{L}
\mathcal{L}(W^T)= \langle y_1^1, w_1(0)\rangle_{H_1,H_{-1}} -
\langle y_1^0, w_1^{\prime}(0)\rangle_{H_2,H_{-2}} +
\langle y_2^1, w_2(0)\rangle_{H,H} - \langle y_2^0, w_2^{\prime}(0)\rangle_{H_1,H_{-1}} \,,
\end{equation}
\noindent where $W=(w_1,w_2,w_1^{\prime},w_2^{\prime})$ and
$\widetilde{W}=(\widetilde{w_1},\widetilde{w_2},\widetilde{w_1}^{\prime},\widetilde{w_2}^{\prime})$ are respectively solutions of \eqref{WT} and
\eqref{tildeW}. Thanks respectively to the admissibility inequality \eqref{directweak1bis} and to the observability inequalities \eqref{directweak2bis}-\eqref{directweak3bis} with $T$ replacing $0$, $\Lambda$ is continuous and coercive on $X_{-1}$ for $T > T_3$. From the usual energy estimates for the time reverse problem for $Z=\mathcal{A}^{-1}W$, and the conservation of $e_0(Z_1)$ through time we have
\begin{equation}\label{ZZTOP}
e_0(Z_1)(0) + e_1(Z_2)(0) \le C( e_0(Z_1)(T) + e_1(Z_2)(T)) \,.
\end{equation}
\noindent This together with  Lemma~\ref{tech1}-(iv), lead to
$$
e_{-1}(W_1)(0) + e_0(W_2)(0) \le C (e_{-1}(W_1)(T) + e_0(W_2)(T)) \,,
$$
\noindent so that $\mathcal{L}$ is continuous
on $X_{-1}$. Hence, 
thanks to Lax-Milgram Lemma, there exists a unique $W^T \in X_{-1}$ such that

\begin{equation}\label{HUM}
\Lambda(W^T,\widetilde{W}^T)=-\mathcal{L}(\widetilde{W}^T) \,, \quad \forall \widetilde{W}^T \in X_{-1} \,.
\end{equation}
\noindent We set $v=B^{\ast}w_2$. Then, thanks to the hidden regularity property
due to  \eqref{directweak1bis}, we have $v \in L^2([0,T];G)$ . Thus, we have by definition of
the solution of \eqref{CTH} by transposition

\begin{multline}\label{TRANS}
\int_0^T \langle v, B^{\ast}\widetilde{w_2}\rangle_G \,dt=
\langle y_1^{\prime}(T), \widetilde{w_1}(T)\rangle_{H_1,H_{-1}} -
\langle y_1(T), \widetilde{w_1}^{\prime}(T)\rangle_{H_2,H_{-2}} +
\langle y_2^{\prime}(T), \widetilde{w_2}(T)\rangle_{H,H} - \\
\langle y_2(T), \widetilde{w_2}^{\prime}(T)\rangle_{H_1,H_{-1}} - 
\mathcal{L}(\widetilde{W}^T) \,, \forall \ 
\widetilde{W}^T \in X_{-1}\,.
\end{multline}
\noindent On the other hand, we have
$$
\int_0^T \langle v, B^{\ast}\widetilde{w_2}\rangle_G \,dt=\Lambda(W^T,\widetilde{W}^T)=-\mathcal{L}(\widetilde{W}^T) \,,
$$
\noindent so that,we deduce from these two relations that $Y(T)=(y_1,y_2,y_1^{\prime},y_2^{\prime})(T)=0$.

\medskip

Assume now that $(ii)$ holds. Let $Y_0 \in X_{1}^{\ast}$.  We consider on
$X_1$ the bilinear form

\begin{equation}\label{Lambdaunb}
\Lambda(U^T,\widetilde{U}^T)= \int_0^T \langle B^{\ast}u_2, B^{\ast}\widetilde{u_2}\rangle_G\,dt \,,
\forall \ U^T, \widetilde{U}^T \in X_1 \,,
\end{equation}
\noindent and the linear form on $X_1$ defined by

\begin{equation}\label{Lunb}
\mathcal{L}(U^T)= \langle y_1^1, u_1(0)\rangle_{H,H} -
\langle y_1^0, u_1^{\prime}(0)\rangle_{H_1,H_{-1}} +
\langle y_2^1, u_2(0)\rangle_{H_{-1},H_1} - \langle y_2^0, u_2^{\prime}(0)\rangle_{H,H} \,,
\forall \  U^T  \in X_1 \,.
\end{equation}
Thanks respectively to the admissibility inequality \eqref{admissineq} and to the observability inequality
\eqref{obsNSH}, $\Lambda$ is continuous and coercive on $X_1$ for $T > T_3$. On the other hand $\mathcal{L}$ is continuous on $X_1$ thanks to \eqref{ZZTOP} with $U$ replacing $Z$. Hence, 
thanks to Lax-Milgram Lemma, there exists a unique $U^T \in X_{1}$ such that

\begin{equation}\label{HUMbd}
\Lambda(U^T,\widetilde{U}^T)=-\mathcal{L}(\widetilde{U}^T) \,, \quad \forall \ \widetilde{U}^T \in X_1 \,.
\end{equation}
\noindent We set $v=B^{\ast}u_2$. We deduce as for the case $(i)$ that $Y(T)=0$.
\end{Proof}
\subsubsection{Further generalizations}
We can change the functional setting in Theorem~\ref{control2} to smoother or weaker control spaces in which null controllability holds, under the appropriate
hypotheses replacing the assumptions $(A2)-(A5)$.  Let $k \in \mathbb{Z}$ be given. We consider the following new set of assumptions.
\begin{equation}\label{hypCk}
(A2)_k
\begin{cases}
\, C_{21}^{\ast} \in \mathcal{L}(H_p) \mbox{ for } p \in \{-k,1-k\} \,,\\
 \, ||A^{k/2}C_{21}A^{-k/2}||=\beta\,,
|A^{k/2}C_{21}A^{-k/2}w|^2 \le \beta \langle A^{k/2}C_{21}A^{-k/2}w\,,w\rangle \quad \forall \ w \in H \,,\\
\exists \alpha>0 \mbox{ such that, }
\alpha\, |\Pi w|^2 \le \langle A^{k/2}C_{21}A^{-k/2}w,w\rangle \quad
\forall \ w \in H \,,
\end{cases}
\end{equation}
\noindent where the operator $\Pi$ satisfies the assumptions $(A3)$.

We make the following assumptions on the observability operator $\mathcal{\mathbf{B}}^{\ast}$.
\begin{equation}\label{admissibilityk}
(A4)_k
\begin{cases}
\mathcal{\mathbf{B}}^{\ast} \in \mathcal{L}(H_{k+2}\times H_k;G),\\
\forall \ T > 0 \ \exists \ C >0 , \mbox{ such that for all } (w^0,w^1) \in H_{k+1} \times H_k \mbox{ and all } f  \in L^2([0,T];H_k)\,,\\
\mbox{the solution } w \mbox{ of }
w'' + A w = f \,, (w,w^{\prime})(0)=  (w^0,w^1)
\mbox{ satisfies }  \\
\int_0^T \| \mathcal{\mathbf{B}}^{\ast}(w,w^{\prime}) \|_G^2 dt \leq C \Big( e_{k+1}(W)(0) + e_{k+1}(W)(T) + \int_0^T e_{k+1}(W)(t) dt + \\
\int_0^T |A^{k/2}f|^2 dt \Big)\,,
\end{cases}
\end{equation}
\noindent where $W=(w,w^{\prime})$.

We further assume the following observability inequality for a single equation
\begin{equation}\label{observabilityBk}
(A5)_k
\begin{cases}
\exists \ T_0>0, \forall \ T>T_0 \,, \exists \ C_1(T)>0 \mbox{ such that }\\
\forall \  (w^0,w^1) \in H_{k+1} \times H_k \,,
\mbox{ the solution } w \mbox{ of }\\
w'' + A w = 0   \,,     (w,w^{\prime})(0)=  (w^0,w^1)
\mbox{ satisfies } \\
\int_0^T \| \mathcal{\mathbf{B}}^{\ast}(w,w^{\prime}) \|_G^2 dt \geq C_1(T)e_{k+1}(W)(0)\,. 
\end{cases}
\end{equation}
Then we can prove the following result.
\begin{Theorem}\label{control2k}
Let $k \in \mathbb{Z}$ be given. We define the sets
$$
X_{-1,k}^{\ast}=H_{2-k}\times H_{1-k}\times H_{1-k} \times H_{-k}\,,
\ X_{1,k}^{\ast}=H_{1-k}\times H_{-k}\times H_{-k} \times H_{-k-1}\,.
$$
\noindent Assume the hypotheses $(A1)$,  $(A2)_k$, $(A3)$ and $(A4)_k-(A5)_k$. We define $T_3>0$ as in Theorem~\ref{NSobs}. Then, we have the following properties.
\begin{itemize}
\item (i) Let $\mathcal{\mathbf{B}}^{\ast}(w,w^{\prime})=B^{\ast}w^{\prime}$ 
where $B \in \mathcal{L}(G,H_k^{\prime})$ is such that $(A4)_k-(A5)_k$ holds. Then, for  all $T >T_3$ ,
and all $Y_0 \in X_{-1,k}^{\ast}$, there exists a control function $v \in L^2((0,T);G)$ such that the solution $Y=(y_1,y_2,y_1^{\prime},y_2^{\prime})$ of \eqref{CTH} satisfies $Y(T)=0$.
\item (ii) Let $\mathcal{\mathbf{B}}^{\ast}(w,w^{\prime})=B^{\ast}w$ where $B \in\mathcal{L}(G, H_{k+2}^{\prime})$ is such that $(A4)_k-(A5)_k$ holds.  
Then, for all $T >T_3$,
and all $Y_0 \in X_{1,k}^{\ast}$, there exists a control function $v \in L^2((0,T);G)$ such that the solution $Y=(y_1,y_2,y_1^{\prime},y_2^{\prime})$ of \eqref{CTH} satisfies $Y(T)=0$.
\end{itemize}
\end{Theorem}
We only sketch the proof since it follows that of the previous section on admissibility and observability for the dual cascade
system, and that of Theorem~\ref{control2} with the appropriate changes of unknowns, of the observability operator and of the coupling term.
\begin{Proof}
We define an operator $\overline{\mathcal{\mathbf{B}}}^{\ast} \in \mathcal{L}(H_2\times H;G)$ as follows.
$$
\overline{\mathcal{\mathbf{B}}}^{\ast}(\overline{u},\overline{v})=\mathcal{\mathbf{B}}^{\ast}(A^{-k/2}\overline{u},A^{-k/2}\overline{v})
\quad \forall \  (\overline{u},\overline{v}) \in H_2 \times H \,.
$$ 

We consider the dual problem
\begin{equation}\label{NSHk}
\begin{cases}
u_1^{\prime\prime} + A u_1 = 0 \,,\\
u_2^{\prime\prime} + A u_2+ C_{21}u_1  = 0 \,,\\
(u_i,u_i^{\prime})(0)=(u_i^0,u_i^1) \in H_{k+2}\times H_k\mbox{ for }
i=1,2 \,.
\end{cases}
\end{equation}
\noindent We introduce the new unknowns $\overline{u}_i=A^{k/2}u_i$ for $i=1,2$. We set $\overline{C}_{21}=A^{k/2}C_{21}A^{-k/2}$.
Then we have $\overline{\mathcal{\mathbf{B}}}^{\ast}\overline{U}_2=\mathcal{\mathbf{B}}^{\ast}U_2$,
where $\overline{U}_2=(\overline{u}_2,\overline{u}_2^{\prime})$. Moreover $\overline{U}=(\overline{u}_1, \overline{u}_2,\overline{u}_1^{\prime},
\overline{u}_2^{\prime})$ solves \eqref{NSH} in $\mathcal{H}$, where $\overline{C}_{21}$ replaces $C_{21}$. By construction
and thanks to $(A4)_k-(A5)_k$, the observability
operator $\overline{\mathcal{\mathbf{B}}}^{\ast}$ satisfies $(A4)-(A5)$ with $f$ replaced by $\overline{f}=A^{k/2}f$
and $W=(w, w^{\prime})$ replaced by $\overline{W}=(A^{k/2}w, A^{k/2}w^{\prime})$ in $(A4)$. Moreover $\overline{C}_{21}$ and $\Pi$ satisfy respectively $(A2)$ and $(A3)$.
Hence we can apply Lemma~\ref{admiss} and Theorem~\ref{NSobs} to $\overline{U}$
with the observability operator $\overline{\mathcal{\mathbf{B}}}^{\ast}$ and the coupling $\overline{C}_{21}$. 

Coming back to the unknown $U$ and
to the observability operator $\mathcal{\mathbf{B}}^{\ast}$, we deduce the following admissibility property
 \begin{equation}\label{admissineqk}
 \int_0^T ||\mathcal{\mathbf{B}}^{\ast}U_2||_G^2 \,dt \le C \Big(e_k(U_1)(0) +e_{k+1}(U_2)(0)\Big) \,,
 \end{equation}
\noindent for all initial data $U^0 \in H_{k+1}^2 \times H_k^2$, and the following observability inequalities
for all $T>T_3$,
and all initial data $U^0 \in H_{k+1}^2 \times H_k^2$,
\begin{equation}\label{obsNSHk}
\begin{cases}
d_1(T) \int_0^T||\mathcal{\mathbf{B}}^{\ast}U_2||_G^2
\ge  e_k(U_1)(0)\,,\\
d_2(T)\int_0^T||\mathcal{\mathbf{B}}^{\ast}U_2||_G^2
\ge  e_{k+1}(U_2)(0) \,, 
\end{cases}
\end{equation}
\noindent where $T_3$, $d_i(T)$ for $i=1,2$ are as in Theorem~\ref{NSobs}. 

This allows us to show that the HUM operator
is well-defined in the case $(i)$ on the set 
$$
X_{-1,k}=H_{k-1}\times H_k\times H_{k-2} \times H_{k-1}\,,
$$
\noindent and in the case $(ii)$ on the set
$$
X_{1,k}=H_{k}\times H_{k+1}\times H_{k-1} \times H_{k}\,,
$$
\noindent By duality, the control problem is well-posed in $X_{-1,k}^{\ast}$ (resp. in $X_{1,k}^{\ast}$) for the case $(i)$ (resp. $(ii)$).
\end{Proof}
\begin{Remark}
It is easy to provide examples of operators $C_{21}$ satisfying $(A2)_k$. It is for instance sufficient to choose $C_{21}=A^{-k/2}DA^{k/2}$ where $D$ satisfies $(A2)$.
\end{Remark}

\begin{Remark}\label{source1}
A further generalization of Theorem~\ref{control2} is  to consider a more general control system of the form \eqref{CTH}, but including additionally source terms $\xi_1$ and $\xi_2$ in appropriate spaces (for well-posedness) on respectively the first and second equation of \eqref{CTH} as considered in subsection 3.3. The above controllability Theorem~\ref{control2} can then be extended with no difficulties to handle these source terms, modifying the HUM method as in the proof of Theorem~\ref{thminsens}.
\end{Remark}
\section{Main applicative results}
We shall give in this section the main consequences on the most well-known examples of applications, namely: wave-type, heat-type and Schr\"odinger cascade coupled systems.

Let $\Omega$ be a bounded open set in $\mathbb{R}^d$ with a sufficiently smooth boundary $\Gamma$. The set $\Omega$ can also be a smooth connected compact Riemannian manifold with or without boundary as in \cite{alaleau11}.  Let $T$ be a given positive time. We recall
the following definition for the Geometric Control Condition of Bardos Lebeau Rauch~\cite{blr92}.
\begin{Definition}
We say that an open subset $\omega$ of $\Omega$ satisfies $(GCC)$ if there exists a time $T>0$ such that every generalized geodesic traveling at speed $1$ in $\Omega$ meets $\omega$ at a time $t<T$. We say that a subset $\Gamma_1$ of the boundary $\Gamma$ satisfies $(GCC)$ if there exists a time $T>0$ such that every generalized geodesic traveling at speed $1$ in $\Omega$ meets $\Gamma_1$ at a time $t<T$ in a non-diffractive point.
\end{Definition}
\nopagebreak[4]
\subsection{Boundary and localized observability/controllability of $2$-coupled cascade wave equations with localized couplings}
\subsubsection{Main observability results for $2$-coupled cascade wave systems}
 We consider the following $2$-coupled cascade system of wave equations
\begin{equation}\label{mixedbinondiagonalH}
\begin{cases}
u_{1,tt} -\Delta  u_1 = 0 \quad \mbox{ in } (0,T)\times \Omega\,,\\
u_{2,tt} -\Delta  u_2+ c_{21}(x)u_1  = 0 \quad \mbox{ in } (0,T)\times \Omega\,,\\
u_i=0 \mbox{ for } i=1, 2 \quad \mbox{ in } (0,T)\times \partial \Omega\,,\\
(u_i,u_{i,t})(0)=(u_i^0,u_i^1) \mbox{ for }
i=1, 2 \quad \mbox{ in } \Omega\,,
\end{cases}
\end{equation}
\noindent where the subscript $t$  denotes the partial derivative  with respect to time $t$. 
We set $H=L^2(\Omega)$, and we consider the operator $A$ defined by $Au=-\Delta u$
for $u \in D(A)=H^2(\Omega) \cap H^1_0(\Omega)$. We set $U^0=
(u_1^0, u_2^0, u_1^1, u_2^1)$ for all the sequel of this section. We make the following assumptions on the coefficient $c_{21}$.
$$
(H1)
\begin{cases}
c_{21} \in W^{1,\infty}(\Omega)\,\\
c_{21} \ge 0 \mbox{ on } \Omega \,,\\
\{c_{21} >0\} \supset \overline{O_2} \mbox{ for some open subsets } O_2 \subset \Omega 
\,.
\end{cases}
$$
\begin{Remark}\label{localizedregion}
The coupling term  $c_{21}$ is
assumed to be strictly positive on $\overline{O_2}$, so that the coupling effect due to this term is effective in a neighborhood of this set. We will say in all the sequel that the subset $\overline{O_2}$ is the region on which the coupling is localized.
\end{Remark}
\noindent We recall that 
\begin{equation}\label{iterateDA}
D(A^{k/2})=\{ u \in H^k(\Omega) \,, u=\Delta u = \ldots \Delta^{[(k-1)/2]}u=0 \mbox{ on } \Gamma\}
\end{equation}
\noindent for all $k \in \{0, 1,\ldots\}$, where $[x]$ stands for the integer part of the real number $x$.

 We do not recall below the admissibility property for the coupled system which allows to show that the observations
of the solution are well-defined (hidden regularity result) in a classical way. This property is given in the abstract Theorem~\ref{NSobs}.
\begin{Theorem}[Observability estimates]\label{ncoupledwavethm}
We assume that  the hypothesis $(H1)$ holds for some open subset $O_2
\subset \Omega$  that satisfies $(GCC)$. Then we have the following results.
\begin{itemize}
\item $(i)$ Locally distributed observation. Let $b_{2}$  be a given function defined on $\Omega$ such that
$$
(H2)
\begin{cases}
b_{2} \in \mathcal{C}(\overline{\Omega}) \,, b_{2} \ge 0 \mbox{ on } \Omega\,,\\
\{b_{2} >0\} \supset \overline{\omega_{2}} \mbox{ for some open subset } \omega_{2} \subset \Omega 
\,,\\
\end{cases}
$$  
\noindent where $\omega_2$ satisfies $(GCC)$. Then there exists $T^{\ast}>0$ such that for
all $T>T^{\ast}$, there exist constants $c_{i,2}(T)>0$, $i=1,2$ such that 
for all $\displaystyle{U^0 \in L^2(\Omega) \times H^1_0(\Omega) \times H^{-1}(\Omega) \times L^2(\Omega)}$,
the following observability inequalities hold

$$
\displaystyle{c_{1,2} (T) ||(u_1^0,u_1^1)||^2_{L^2(\Omega)\times H^{-1}(\Omega))} \le \int_0^T \int_{\Omega}b_2 |u_{2,t}|^2 \,dx\, dt  }\,,  
$$
$$
\displaystyle{c_{2,2} (T) ||(u_2^0,u_2^1)||^2_{H^1_0(\Omega) \times L^2(\Omega)} \le \int_0^T \int_{\Omega}b_2 |u_{2,t}|^2 \,dx\, dt  }\,. 
$$
\item $(ii)$ Boundary observability. Let $b_{2}$  be a given function defined on 
$\Gamma$  such that
$$
(H3)
\begin{cases}
b_{2} \in \mathcal{C}(\overline{\Gamma}) \,, b_{2} \ge 0 \mbox{ on } \Gamma \mbox{ for } \,,\\
\{b_{2} >0\} \supset \overline{\Gamma_{2}} \mbox{ for some subset } \Gamma_{2} \subset \Gamma
\,,\\
\end{cases}
$$  
\noindent where $\Gamma_{2}$ satisfies $(GCC)$. Then there exists $T^{\ast}>0$ such that for
all $T>T^{\ast}$, there exist constants $c_{i,2}(T)>0$, $i=1,2$ such that 
for all $\displaystyle{U^0 \in L^2(\Omega) \times H^1_0(\Omega) \times H^{-1}(\Omega) \times L^2(\Omega)}$ the following observability inequalities hold

$$
\displaystyle{c_{1,2} (T) ||(u_1^0,u_1^1)||^2_{L^2(\Omega)\times H^{-1}(\Omega))} \le \int_0^T \int_{\Gamma} b_2\Big|\frac{\partial u_{2}}{\partial \nu}\Big|^2 \,d \sigma\, dt } \,, 
$$
$$
\displaystyle{c_{2,2} (T) ||(u_2^0,u_2^1)||^2_{H^1_0(\Omega) \times L^2(\Omega)} \le \int_0^T \int_{\Gamma} b_2\Big|\frac{\partial u_{2}}{\partial \nu}\Big|^2 \,d \sigma\, dt }\,. 
$$
\end{itemize}
\end{Theorem}
\begin{Remark}
As in the Remark~\ref{localizedregion}, the subsets $\omega_{2}$ (resp. $\Gamma_{2}$) of $\omega$ (resp. $\Gamma$) are the regions (indeed in a neighborhood of them) on which the observations are localized. Hence, we prove that $2$-coupled cascade systems
with a coupling term localized on a subregion $O_2 \subset \Omega$ with a single observation either  locally distributed on $\omega_2 \subset \Omega$ or distributed on a part $\Gamma_2 \subset \Gamma$ of the boundary, is observable under the geometric condition that both $O_2$ and $\omega_2$
(resp. $O_2$ and $\Gamma_2$) satisfy (GCC). Hence this covers many situations for which the intersection
$\overline{O_2} \cap \overline{\omega_2}=\emptyset$ (resp. $\overline{O_2} \cap \overline{\Gamma_2}=\emptyset$).  As shown in Theorem~\ref{2NEC} --which gives a necessary and sufficient condition for the observability of $2$-coupled system-- this condition is sharp, and even necessary when $\partial \Omega$ has no contact of infinite order with its tangent, or is analytic. Moreover, we can exhibit many situations for which this sharp condition holds with
$\overline{O_2}\cap \overline{\omega_2}=\emptyset$ (resp. $\overline{O_2}\cap \overline{\Gamma_2}=\emptyset$).
\end{Remark}
\subsubsection{Main controllability results for $2$-coupled cascade wave systems}
We now consider the $2$-coupled control cascade system subjected to a single either locally distributed or a boundary control.
We shall start with the case of a locally distributed control.
\begin{equation}\label{nmixedwavecont}
\begin{cases}
y_{1,tt} -\Delta y_1 +c_{21}y_2= 0 \mbox{ in } (0,T) \times \Omega \,,\\
y_{2,tt} -\Delta  y_{2}=b_{2} v_{2} \mbox{ in } (0,T) \times \Omega\,\\
y_i=0 \mbox{ in } (0,T) \times \Gamma  \mbox{ for }
i=1, 2 \,, \\
(y_i,y_{i,t})(0)=(y_i^0,y_i^1) \mbox{ for }
i=1, 2 \,.
\end{cases}
\end{equation}
\noindent  We set
$Y_0=(y_1^0, y_{2}^0, y_1^1, y_{2}^1)$. Then we have the following exact controllability result.

\begin{Theorem}\label{contwave}
We assume that the coefficient $c_{21}$ satisfy the hypothesis $(H1)$ for some open subsets $O_2
\subset \Omega$ that satisfies $(GCC)$. We also assume that the coefficients $b_{2}$ and the subset $\omega_{2}$ satisfies $(H2)$ where the subsets $\omega_{2}$  satisfies $(GCC)$. Then for all
$T>T^{\ast}$, and all $Y_0 \in \displaystyle{\Pi_{i=1}^{2} D(A^{(3-i)/2}) \times \Pi_{i=1}^{2} D(A^{(2-i)/2})}$, there exists a control function $v_{2} \in L^2((0,T);L^2(\Omega))$
such that the solution $Y=(y_1, y_{2},y_1^{\prime},y_{2}^{\prime})$ of \eqref{nmixedwavecont} with initial data $Y_0$ satisfies $Y(T)=0$.
\end{Theorem}

We now consider the following $2$-coupled control cascade system with a boundary control.
\begin{equation}\label{nmixedwavecontbd}
\begin{cases}
y_{1,tt} -\Delta y_1 +c_{21}y_2= 0  \mbox{ in } (0,T) \times \Omega\,,\\
y_{2,tt} -\Delta  y_{2}=0 \mbox{ in } (0,T) \times \Omega\,\\
y_1=0 \mbox{ in } (0,T) \times \Gamma \,, y_{2}=b_{2} v_{2}  \mbox{ in } (0,T) \times \Gamma  \,, \\
(y_i,y_{i,t})(0)=(y_i^0,y_i^1) \mbox{ for }
i=1, \ldots, 2 \,.
\end{cases}
\end{equation}
\noindent Then we have the following exact controllability result.

\begin{Theorem}\label{contwavebd}
We assume that the coefficient $c_{21}$ satisfies the hypothesis $(H1)$ for some open subset $O_2
\subset \Omega$ that satisfies $(GCC)$. We also assume that the coefficient $b_{2}$ and the subset $\Gamma_{2}$ satisfy $(H3)$ where the subset $\Gamma_{2}$  satisfies $(GCC)$. Then for all
$T>T^{\ast}$, and all $Y_0 \in \displaystyle{\Pi_{i=1}^{2} D(A^{(2-i)/2}) \times \Pi_{i=1}^{2} 
D(A^{(1-i)/2})}$,
 there exists a control function $v_{2} \in L^2((0,T);L^2(\Gamma))$
such that the solution $Y=(y_1,y_{2},y_1^{\prime}, y_{2}^{\prime})$ of \eqref{nmixedwavecontbd} with initial data $Y_0$ satisfies $Y(T)=0$.
\end{Theorem}
\subsection{Boundary and localized controllability of $2$-coupled cascade heat and Schr\"odinger equations with localized couplings}
We now consider the following $2$-coupled locally control cascade heat-type system. We shall first consider the case of a locally distributed control.
\begin{equation}\label{nmixediffcont}
\begin{cases}
e^{i \theta} y_{1,t} -\Delta y_1 +c_{21}y_2= 0 \mbox{ in } (0,T) \times \Omega \,,\\
e^{i \theta}y_{2,t} -\Delta  y_{2}=b_{2} v_{2} \mbox{ in } (0,T) \times \Omega\,\\
y_i=0 \mbox{ in } (0,T) \times \Gamma  \mbox{ for }
i=1, 2 \,, \\
y_i(0)=y_i^0 \mbox{ in } \Omega \mbox{ for }
i=1, 2 \,.
\end{cases}
\end{equation}
\noindent We set
$Y_0=(y_1^0, y_{2}^0)$. We recover $2$-coupled heat (resp. Schr\"odinger) cascade systems when $\theta=0$ (resp. $\theta=\pm \pi/2$) and diffusive coupled cascade systems when
$\theta \in (-\pi/2, \pi/2)$. 
Then we have the following exact controllability result.

\begin{Corollary}\label{contheat}
We assume that the coefficient $c_{21}$ satisfies the hypothesis $(H1)$ for some open subset $O_2
\subset \Omega$ that satisfies $(GCC)$. We also assume that the coefficient $b_{2}$ and the subset $\omega_{2}$ satisfy $(H2)$ where the subset $\omega_{2}$  satisfies $(GCC)$. Then, 
the following properties hold
\begin{itemize}
\item (i) The case $\theta \in (-\pi/2,\pi/2)$ (Heat type systems). We have  for all
$T>0$, and all $Y_0 \in ( L^2(\Omega))^{2}$, there exist a control function $v_{2} \in L^2((0,T);L^2(\Omega))$
such that the solution $Y=(y_1,y_{2})$ of \eqref{nmixediffcont} with initial data $Y_0$ satisfies $Y(T)=0$.
\item (ii) The case $\theta=\pm \pi/2$ (Schr\"odinger systems). We have for all $T>0$ and all
$Y_0 \in \displaystyle{H^1_0(\Omega) \times L^2(\Omega)}$, there exist a control function $v_{2} \in L^2((0,T);L^2(\Omega))$
such that the solution $Y=(y_1, y_{2})$ of \eqref{nmixediffcont} with initial data $Y_0$ satisfies $Y(T)=0$.
\end{itemize}
\end{Corollary}

We now consider the following $2$-coupled control cascade heat-type system with a boundary control.
\begin{equation}\label{nmixediffcontbd}
\begin{cases}
e^{i \theta}y_{1,t} -\Delta y_1 +c_{21}y_2= 0  \mbox{ in } (0,T) \times \Omega\,,\\
e^{i \theta}y_{2,t} -\Delta  y_{2}=0 \mbox{ in } (0,T) \times \Omega\,\\
y_1=0 \mbox{ in } (0,T) \times \Gamma \,, y_{2}=b_{2} v_{2}  \mbox{ in } (0,T) \times \Gamma   \,, \\
y_i(0)=y_i^0 \mbox{ in } \Omega \mbox{ for }
i=1, 2 \,.
\end{cases}
\end{equation}
\noindent Then we have the following exact controllability result.

\begin{Corollary}\label{contdiffbd}
We assume that the coefficients $c_{21}$ satisfies the hypothesis $(H1)$ for some open subset $O_2
\subset \Omega$ that satisfies $(GCC)$. We also assume that the coefficients $b_{2}$ and the subsets $\Gamma_{2}$ satisfy $(H3)$ where the subset $\Gamma_{2}$  satisfies $(GCC)$. Then 
we have
\begin{itemize}
\item (i) The case $\theta \in (-\pi/2,\pi/2)$ (Heat type systems). We have  for all
$T>0$, and all $Y_0 \in ( H^{-1}(\Omega))^{2}$, there exist a control function $v_{2} \in L^2((0,T);L^2(\Gamma))$
such that the solution $Y=(y_1,y_{2})$ of \eqref{nmixediffcontbd} with initial data $Y_0$ satisfies $Y(T)=0$.
\item (ii) The case $\theta=\pm \pi/2$ (Schr\"odinger systems). We have for all $T>0$ and all
$Y_0 \in \displaystyle{L^2(\Omega) \times H^{-1}(\Omega)}$, there exist a control function $v_{2} \in L^2((0,T);L^2(\Omega))$
 such that the solution $Y=(y_1, y_{2})$ of \eqref{nmixediffcontbd} with initial data $Y_0$ satisfies $Y(T)=0$.
\end{itemize}
\end{Corollary}
\begin{Remark}
It should be noted that $(GCC)$ is not a natural condition for the scalar heat equation for which null-controllability holds for arbitrary nonempty control region. This property probably still holds for coupled cascade system. In fact, when $O_2 \cap \omega_2 \neq \emptyset$, condition GCC is not needed (see~\cite{GBdT10}). Hence the results presented here are not optimal as far as geometric conditions are concerned. They are optimal for one-dimensional domains. Since Jaffard's results~\cite{jaffard} on internal controllability/observability of Schr\"odinger equation in a rectangle, it is also known that $(GCC)$ is not a necessary condition for Schr\"odinger equation (see also~\cite{tenentucs}). However the results presented above are valid in a multidimensional setting, for boundary as well as locally distributed control and for geometric situations for which the control region does not meet the coupling region.
\end{Remark}
\begin{Remark}
The above observability and controllability results also hold true if the operator $-\Delta$ is replaced by a more general uniformly elliptic operator, or by higher order operators such as the bilaplacian for instance. 
\end{Remark}

\subsection{Insensitizing controls for the wave equation}
In this part we give new results on the existence of insensitizing controls for the wave equation, which can be either locally distributed controls or boundary controls. This problem reads as follows. We consider the scalar wave equation with a locally distributed control $v$ (the localization of the control depends on the the support of
the coefficient function $b$):
\begin{equation}\label{insen1}
\begin{cases}
y_{tt} - \Delta y = \xi + bv \quad \mbox{ in } (0,T) \times \Omega\,,\\
y=0 \quad \mbox{ in } (0,T) \times \Gamma\,,\\
y(0,\cdot)=y^0 + \tau_0 z^0 \mbox{ in } \Omega \,, y_t(0,\cdot)=y^1 + \tau_1 z^1 \mbox{ in } \Omega \,,
\end{cases}
\end{equation}
\noindent and the scalar wave equation with a boundary  control $v$ (the localization of the control depends on the the support of
the coefficient function $b$ in $\Gamma$):
\begin{equation}\label{insen2}
\begin{cases}
y_{tt} - \Delta y = \xi  \quad \mbox{ in } (0,T) \times \Omega\,,\\
y=b v \quad \mbox{ in } (0,T) \times \Gamma\,,\\
y(0,\cdot)=y^0 + \tau_0 z^0 \mbox{ in } \Omega \,, y_t(0,\cdot)=y^1 + \tau_1 z^1 \mbox{ in } \Omega \,,
\end{cases}
\end{equation}
\noindent where for both cases, the inhomogeneity $\xi \in L^2((0,T)\times \Omega)$, the initial data $(y^0,y^1)$ are given known functions in $H^1_0(\Omega)\times L^2(\Omega)$ or in $L^2(\Omega)\times H^{-1}(\Omega)$, whereas the
perturbations $z^0$, $z^1$ are supposed to be unknown of
norm $1$ in the appropriate spaces. Here the real numbers $\tau_0, \tau_1$ are assumed to be small  and to measure the amplitudes of the unknown perturbations of the initial data. We associate to the solutions $y$ of \eqref{insen1} (resp. to \eqref{insen2}) the functional defined by
\begin{equation}\label{func}
\Phi(y, \tau_0, \tau_1)= \frac{1}{2} \int_0^T\int_{\Omega}c  y^2 \,dx dt \,,
\end{equation}
\noindent where $c$ is a function which has a localized support in the neighborhood of a subset $O$ which is a given subset of $\Omega$. The functional $\Phi$ consists in a ``weighted" observation of the solution on the set $O$ during a length of time $T$. It is an observation when the coefficient $c=\mathds{1}_O$. The control $bv$ is said to insensitize the
``weighted" observation $\Phi$ if for all $(z^0, z^1)$ the corresponding solution $y$ of \eqref{insen1} (resp. \eqref{insen2}) satisfies
\begin{equation}\label{insen3}
\frac{\partial \Phi}{\partial \tau_0}(y,0,0)=\frac{\partial \Phi}{\partial \tau_1}(y,0,0)=0 \,.
\end{equation}
\noindent  We shall see that \eqref{insen3} is equivalent to solve an appropriate exact controllability result for a cascade system, namely for \eqref{cascade1}
(resp. \eqref{cascade2}) for the locally distributed (resp. boundary) insensitizing control problem.
We refer to J.-L. Lions~\cite{lions89} for the introduction of this notion and its reduction to a controllability problem of an associated cascade system. This notion was further weakened to a notion of $\varepsilon$-insensitizing control for parabolic equations by Bodart and Fabre~\cite{bodart-fabre95}. The first results concerning the existence of boundary insensitizing (resp.
$\varepsilon$- localized insensitizing) controls for the wave equation have been obtained by D\'ager~\cite{Dager06} in the case of the one-dimensional wave equation. They have been extended by Tebou~\cite{tebou08} to the multidimensional wave equation in the case of locally distributed control, and for coupling regions which meet the control region. The case of localized insensitizing control for a semilinear wave equation with a different observation (namely the trace of the normal derivative of the solution) has been considered in \cite{tebou2011}. Positive results have been obtained in  \cite{tebou2011} under the strong geometric condition that the boundary observation region satisfies the usual multiplier geometric condition and that the localized control region contains a neighbourhood of the observation region.

Here we are interested in general results for \eqref{insen1} and \eqref{insen2} in sharp geometric situations, in particular we want positive results in situations for which $supp\{c\} \cap supp\{b\} =\emptyset$.

\begin{Theorem}\label{thminsens}
Assume that $c$ satisfies
\begin{equation}\label{cinsen}
\begin{cases}
c \in W^{1,\infty}(\Omega)\,, c \ge 0 \mbox{ on } \Omega\,,\\
\{c >0\} \supset \overline{O} \mbox{ for some open subset } O \subset \Omega  \,.
\end{cases}
\end{equation}
\noindent  We have the following properties
\begin{itemize}
\item Locally distributed control. Let $b$ be given such that
\begin{equation}\label{binsen}
\begin{cases}
b \in \mathcal{C}(\overline{\Omega})\,, b \ge 0 \mbox{ on } \Omega\,,\\
\{b >0\} \supset \overline{\omega} \mbox{ for some open subset } \omega \subset \Omega  \,.
\end{cases}
\end{equation}
\noindent We assume that $O$ and $\omega$ satisfy $(GCC)$. Then for any given $\xi \in L^2((0,T);L^2(\Omega))$ and $(y^0,y^1) \in H^1_0(\Omega)\times L^2(\Omega)$, there exists
an exact control  $v \in L^2((0,T);L^2(\Omega))$ that  drives back the solution of the scalar wave equation
$$
\begin{cases}
y_{2,tt} -\Delta y_2 =\xi + b v \quad \mbox{ in } (0,T) \times \Omega \,,\\
y_2=0 \quad \mbox{ in } (0,T) \times \Gamma \,,\\
(y_2, y_{2,t})_{| t=0}=(y^0,y^1) \quad \mbox{in } \Omega \,,
\end{cases}
$$
to equilibrium, i.e. $y_2(T)=y_{2,t}(T)=0$ and insensitizes $\Phi$ along the solutions of \eqref{insen1},
i.e. $v$ is such that  \eqref{insen3} holds for any $(z^0,z^1) \in
H^1_0(\Omega) \times L^2(\Omega)$.
\item Boundary control. Let $b$ be given such that
\begin{equation}\label{binsenbd}
\begin{cases}
b \in \mathcal{C}(\overline{\Gamma})\,, b \ge 0 \mbox{ on } \Gamma\,,\\
\{b >0\} \supset \overline{\Gamma_1} \mbox{ for some subset } \Gamma_1 \subset \Gamma  \,.
\end{cases}
\end{equation}
\noindent We assume that $O$ and $\Gamma_1$ satisfy $(GCC)$. Then for any given $\xi \in L^2((0,T);L^2(\Omega))$ and $(y^0,y^1) \in L^2(\Omega) \times H^{-1}(\Omega)$, there exists
an exact control control  $v \in L^2((0,T);L^2(\Gamma))$ that drives back the solution of the scalar wave equation
$$
\begin{cases}
y_{2,tt} -\Delta y_2 =\xi \quad \mbox{ in } (0,T) \times \Omega \,,\\
y_2=bv \quad \mbox{ in } (0,T) \times \Gamma \,,\\
(y_2, y_{2,t})_{| t=0}=(y^0,y^1) \quad \mbox{in } \Omega \,,
\end{cases}
$$
to equilibrium, i.e. $y_2(T)=y_{2,t}(T)=0$ and insensitizes $\Phi$ along the solutions of \eqref{insen2},
i.e. $v$ is such that  \eqref{insen3} holds for any $(z^0,z^1) \in L^2(\Omega) \times H^{-1}(\Omega)$.
\item Assume that $\partial \Omega$ has no contact of infinite order with its tangent, or is analytic. Then the above  condition on $O$ and $\omega$ (resp. on $O$ and $\Gamma_1$) for the case of locally (resp. boundary) distributed control are necessary.
\end{itemize}
\end{Theorem}
\section{Proofs of the main applicative results}
\subsection{Proofs of the results for coupled cascade wave systems}
We begin with the proof of Theorem~\ref{ncoupledwavethm}.
\begin{Proof} of Theorem~\ref{ncoupledwavethm}.
This is an application of the abstract Theorem~\ref{NSobs}. Here $H=L^2(\Omega)$, and
$A$ is given by $A=-\Delta$ and $D(A)=H^2(\Omega) \cap H^1_0(\Omega)$. The sets $H_k=
D(A^{k/2})$ are given by \eqref{iterateDA}. The operator $C_{21}$ is the multiplication operator 
in $L^2(\Omega)$ by the corresponding function $c_{21}$ and is therefore bounded and self-adjoint.Thanks to the smoothness assumptions on the  coefficients $c_{21}$, $C_{21} \in \mathcal{L}(H_k)$ for all
$k=0,1$. Thanks to $(H1)$, the assumption $(A2)$ holds with $\Pi=
\mathds{1}_{O_2}$ and $\alpha=\inf_{O_2}(c_{2 1})>0$. Morever since the set $O_2$ satisfies
$(GCC)$,  $(A3)$ also holds. We shall check the assumptions on the observability operator.
 First case: Locally distributed observation. Then we have $G=G_{2}=L^2(\Omega)$,
and $B^{\ast}=B_{2}^{\ast}$ is the multiplication in $L^2(\Omega)$ by the bounded function $b_{2}$. Therefore $B_{2}^{\ast}$ is a bounded symmetric operator in $H$, so that
$(A4)$ holds with $B^{\ast}=B_{2}^{\ast}$. Thanks to the assumptions
on $b_{2}$ in the case of locally distributed observation and since $\omega_{2}$ satisfies $(GCC)$, we deduce
that $(A5)$ also holds using \cite{blr92}.

Second case: Boundary observability. Then we have $G=G_{2}=L^2(\Gamma)$,
and $B^{\ast}=B_{2}^{\ast} \in \mathcal{L}(H^2(\Omega) \cap H^1_0(\Omega);L^2(\Gamma))$ is defined as
$$
B_{2}^{\ast}u=b_{2} \frac{ \partial u}{\partial \nu} \quad u \in H^2(\Omega) \cap H^1_0(\Omega) \,.
$$
\noindent Thanks to the well-known hidden regularity result of~\cite{lions}, $B^{\ast}=B_{2}^{\ast}$ satisfies $(A4)$. Thanks to the assumptions
on $b_{2}$ in the case of boundary observability and since $\Gamma_{2}$ satisfies $(GCC)$, we deduce
that $(A5)$ also holds using \cite{blr92}. Hence we can apply Theorem~\ref{NSobs}, which gives the desired result.
\end{Proof}
We now shall prove the controllability results for coupled cascade wave systems.
\begin{Proof} of Theorem~\ref{contwave}.
Thanks to our hypotheses and thanks to the above proof the assumption of
Theorem~\ref{NSobs} are satisfied. Hence we apply the part $(i)$ of Theorem~\ref{control2}. This gives the desired result.
\end{Proof}
\begin{Proof} of Theorem~\ref{contwavebd}.
Thanks to our hypotheses and thanks to the above proof, the assumption of
Theorem~\ref{NSobs} are satisfied. Hence we apply the part $(ii)$ of Theorem~\ref{control2}. This gives the desired result.
\end{Proof}
\subsection{Proofs of the results for coupled cascade heat and Schr\"odinger equations}
We start with the proof of Corollary~\ref{contheat} using the transmutation method.
\begin{Proof} of Corollary~\ref{contheat}.
We proceed as in~\cite{alaleau11}. 
Proof of $(i)$. We first apply no control on the time interval $(0,T/2)$, so that by
the smoothing effect of the heat equation, the initial data $Y^0\in (L^2(\Omega))^{2}$ is driven to
$Y_{|t=T/2} \in  \displaystyle{\Pi_{i=1}^{2} D(A^{(3-i)/2})}$. We then combine Theorem~\ref{contwave} together with the transmutation result  given by Miller in~\cite{miller04} to prove that there exists a control $v_{2}$ such that $Y(T)=0$.

Proof of $(ii)$. It is similar to the case $(i)$ except that we work directly on the time interval $(0,T)$ since no smoothing effect holds in the case of the Schr\"odinger equation. Combining Theorem~\ref{contwave} together with the transmutation result given by Miller in~\cite{miller05} (see Theorem 3.1), we conclude the proof.
\end{Proof}
\begin{Proof} of Corollary~\ref{contdiffbd}.
We proceed as in~\cite{alaleau11}. 
Proof of $(i)$. We first apply no control on the time interval $(0,T/2)$, so that by
the smoothing effect of the heat equation, the initial data $Y^0 \in (H^{-1}(\Omega))^{2}$ is driven to
$Y_{|t=T/2} \in  \displaystyle{\Pi_{i=1}^{2} D(A^{(2-i)/2})}$. We then combine Theorem~\ref{contwavebd} together with the transmutation result given by Miller in~\cite{miller} (see Theorem 3.4) to get the desired result.

Proof of $(ii)$. It is similar to the case $(i)$ except that we work directly on the time interval $(0,T)$ since no smoothing effect holds in the case of the Schr\"odinger equation. Combining Theorem~\ref{contwavebd} together with the transmutation result (Theorem 3.1) given by Miller in~\cite{miller05} (see Theorem 3.1), we conclude the proof.
\end{Proof}

\subsection{Proofs of the results on insensitizing controls for the wave equation}
The proof of Theorem~\ref{thminsens} relies on the appropriate admissibility and observability inequalities for $2$-coupled cascade systems given in the Lemmata~\ref{admiss}, \ref{obsdir} and in Theorem~\ref{NSobs}. Namely it relies on the  following classical proposition (established in a different form in~\cite{Dager06} for the one-dimensional case) linking the existence of insensitizing controls to that of certain properties for an associated cascade system.
\begin{Proposition}\label{insen4}
Assume $c \in W^{1,\infty}(\Omega)$, and let $(y^0,y^1) \in H^1_0(\Omega)\times L^2(\Omega)$ and $\xi \in L^2((0,T);L^2(\Omega))$ be given. We have the following properties
\begin{itemize}
\item Locally distributed control. Let $b$ be given in $\mathcal{C}(\overline{\Omega})$, and locally supported
in $\overline{\Omega}$. If $v \in L^2((0,T);L^2(\Omega))$ is such that
the solution of the $v$-controlled cascade system
\begin{equation}\label{cascade1}
\begin{cases}
y_{1,tt} -\Delta y_1 + c(x) y_2=0 \quad \mbox{ in } (0,T) \times \Omega \,,\\
y_{2,tt} -\Delta y_2 =\xi + b v \quad \mbox{ in } (0,T) \times \Omega \,,\\
y_1=y_2=0 \quad \mbox{ in } (0,T) \times \Gamma \,,\\
Y^0=(y_1,y_2,y_{1,t}, y_{2,t})_{| t=0}=(0,y^0,0,y^1) \quad \mbox{in } \Omega \,,
\end{cases}
\end{equation}
\noindent satisfies 
\begin{equation}\label{insen0}
y_1(T)=y_{1,t}(T)=0 \mbox{ in } \Omega \,,
\end{equation}
\noindent then $v$ insensitizes $\Phi$ along the solutions of \eqref{insen1}. In particular if the control $v$ can be chosen such that the controlled cascade system \eqref{cascade1} satisfies $Y(T)=(y_1(T),y_2(T),y_{1,t}(T),y_{2,t}(T))=(0,0,0,0)$, then $v$ is an exact control for
the scalar wave equation
\begin{equation}\label{cascadesing1}
\begin{cases}
y_{2,tt} -\Delta y_2 =\xi + b v \quad \mbox{ in } (0,T) \times \Omega \,,\\
y_2=0 \quad \mbox{ in } (0,T) \times \Gamma \,,\\
(y_2, y_{2,t})_{| t=0}=(y^0,y^1) \quad \mbox{in } \Omega \,,
\end{cases}
\end{equation}
\noindent which both drives the solution of the scalar wave equation \eqref{cascadesing1} to equilibrium and insensitizes $\Phi$ along the solutions of \eqref{insen1}.

Conversely, if $v$ insensitizes $\Phi$ along the solutions of \eqref{insen1}, then the solution of \eqref{cascade1} satisfies \eqref{insen0}.
\item Boundary control. Let $b$ be given in $\mathcal{C}(\overline{\Gamma})$, and locally supported
in $\Gamma$. If $v \in L^2((0,T);L^2(\Gamma))$ is such that
the solution of the $v$-controlled cascade system
\begin{equation}\label{cascade2}
\begin{cases}
y_{1,tt} -\Delta y_1 + c(x) y_2=0 \quad \mbox{ in } (0,T) \times \Omega \,,\\
y_{2,tt} -\Delta y_2 =\xi \quad \mbox{ in } (0,T) \times \Omega \,,\\
y_1=0\,, y_2=bv \quad \mbox{ in } (0,T) \times \Gamma \,,\\
Y^0=(y_1,y_2,y_{1,t}, y_{2,t})_{| t=0}=(0,y^0,0,y^1) \quad \mbox{in } \Omega \,,
\end{cases}
\end{equation}
\noindent satisfies \eqref{insen0},
\end{itemize}
\noindent then,
$v$ insensitizes $\Phi$ along the solutions of \eqref{insen2}. In particular if the control $v$ can be chosen such that the controlled cascade system \eqref{cascade2} satisfies $Y(T)=(y_1(T),y_2(T),\linebreak y_{1,t}(T),y_{2,t}(T))=(0,0,0,0)$, then $v$ is an exact control for
the scalar wave equation
\begin{equation}\label{cascadesing2}
\begin{cases}
y_{2,tt} -\Delta y_2 =\xi  \quad \mbox{ in } (0,T) \times \Omega \,,\\
y_2=bv \quad \mbox{ in } (0,T) \times \Gamma \,,\\
(y_2, y_{2,t})_{| t=0}=(y^0,y^1) \quad \mbox{in } \Omega \,,
\end{cases}
\end{equation}
\noindent which both drives the solution of the scalar wave equation \eqref{cascadesing2} to equilibrium and insensitizes $\Phi$ along the solutions of \eqref{insen2}.

Conversely if $v$ insensitizes $\Phi$ along the solutions of \eqref{insen2}, then the solution of \eqref{cascade2} satisfies \eqref{insen0}.
\end{Proposition}
\begin{Proof} 
We first consider the case of locally distributed control. Thanks to our assumptions on $b, \xi, y^0, y^1$ the equation in $y_2$ has a unique
solution $y_2 \in \mathcal{C}([0,T]; H^1_0(\Omega))\cap \mathcal{C}^1([0,T]; L^2(\Omega))$. Thus in a similar
way the equation in $y_1$ has a unique solution $y_1 \in \mathcal{C}([0,T]; H^2(\Omega) \cap H^1_0(\Omega))\cap \mathcal{C}^1([0,T];H^1_0(\Omega))$. Assume that the solution of \eqref{cascade1}
satisfies \eqref{insen0}. We have
\begin{equation}\label{cas3}
\frac{\partial \Phi}{\partial \tau_0}(y,0,0)=\int_0^T\int_{\Omega}cy_2 \hat{w} \,dx\,dt \,,
\end{equation}
\noindent and
\begin{equation}\label{cas4}
\frac{\partial \Phi}{\partial \tau_1}(y,0,0)=\int_0^T\int_{\Omega}c y_2 \hat{z} \,dx\,dt \,,
\end{equation}
\noindent for all $(z^0,z^1) \in H^1_0(\Omega) \times L^2(\Omega)$ (of norm $1$ in these spaces), where $\hat{w}, \hat{z}$ are the respective solutions of
\begin{equation}\label{wcas}
\begin{cases}
\hat{w}_{tt} - \Delta \hat{w} = 0 \quad \mbox{ in } (0,T) \times \Omega\,,\\
\hat{w}=0 \quad \mbox{ in } (0,T) \times \Gamma\,,\\
\hat{w}(0,\cdot)=z^0  \mbox{ in } \Omega \,, \hat{w}_t(0,\cdot)=0  \mbox{ in } \Omega \,,
\end{cases}
\end{equation}
\begin{equation}\label{zcas}
\begin{cases}
\hat{z}_{tt} - \Delta \hat{z} = 0 \quad \mbox{ in } (0,T) \times \Omega\,,\\
\hat{z}=0 \quad \mbox{ in } (0,T) \times \Gamma\,,\\
\hat{z}(0,\cdot)=0  \mbox{ in } \Omega \,, \hat{z}_t(0,\cdot)=z^1  \mbox{ in } \Omega \,.
\end{cases}
\end{equation}
\noindent We set
$$
w(t,\cdot)=\hat{w}(T-t,\cdot) \,, z(t,\cdot)=\hat{z}(T-t,\cdot) \quad \forall \ t \in [0,T] \,.
$$
\noindent Then $w$ and $z$ solves respectively
\begin{equation}\label{wcasrev}
\begin{cases}
w_{tt} - \Delta w= 0 \quad \mbox{ in } (0,T) \times \Omega\,,\\
w=0 \quad \mbox{ in } (0,T) \times \Gamma\,,\\
w(T,\cdot)=z^0  \mbox{ in } \Omega \,, w_t(T,\cdot)=0  \mbox{ in } \Omega \,,
\end{cases}
\end{equation}
\begin{equation}\label{zcasrev}
\begin{cases}
z_{tt} - \Delta z = 0 \quad \mbox{ in } (0,T) \times \Omega\,,\\
z=0 \quad \mbox{ in } (0,T) \times \Gamma\,,\\
z(T,\cdot)=0  \mbox{ in } \Omega \,, z_t(T,\cdot)=-z^1  \mbox{ in } \Omega \,.
\end{cases}
\end{equation}
\noindent We set $\hat{y}_i(t)=y_i(T-t)$ for $t \in [0,T]$ and $i=1,2$. Then, by definition of the solution by transposition (see \eqref{TRANSy1} in Remark~\ref{rktrans})
and thanks to the relations $\hat{y}_1(T)=\hat{y}_{1,t}(T)=0$, we have  
$$
-\int_0^T\langle \hat{y}_2 \,, cw\rangle_{H_1,H_{-1}}\,dt=-
\langle \hat{y}_{1,t}(0), w(0)\rangle_{H_1,H_{-1}} +
\langle \hat{y}_1(0), w_t(0)\rangle_{H_2,H_{-2}}\,.
$$
\noindent Since $\hat{y}_1(0)=\hat{y}_{1,t}(0)=0$ and $w(t) \in L^2(\Omega)$ for all $t \in [0,T]$, we deduce that  
$$
\int_0^T\int_{\Omega}c \hat{y}_2 w\,dx\,dt=\int_0^T\int_{\Omega}c y_2 \hat{w}\,dx\,dt=0 \quad \forall \ z^0 \in H^1_0(\Omega) \,.
$$
In a similar way, replacing $w$ by $z$ in the above relation, we deduce that
$$
\int_0^T\int_{\Omega}cy_2\hat{z}\,dx \,dt=0 \quad \forall \ z^1 \in L^2(\Omega)\,.
$$
\noindent Hence $v$ insensitizes $\Phi$ along the solutions of \eqref{insen1}.

\noindent Conversely assume that $v$ insensitizes $\Phi$ along the solutions of \eqref{insen1}. We denote
by $y_2$ the solution of
$$
\begin{cases}
y_{2,tt} -\Delta y_2 =\xi +bv \quad \mbox{ in } (0,T) \times \Omega \,,\\
 y_2=0 \quad \mbox{ in } (0,T) \times \Gamma \,,\\
(y_2,y_{2,t})_{| t=0}=(y^0,y^1) \quad \mbox{in } \Omega \,.
\end{cases}
$$
\noindent Then 
$$
\int_0^T\int_{\Omega}cy_2 \hat{w} \,dx\,dt=\int_0^T\int_{\Omega}c y_2 \hat{z} \,dx\,dt=0 \,,
$$
\noindent for all $(z^0,z^1) \in H^1_0(\Omega) \times L^2(\Omega)$ (of norm $1$ in these spaces), where $\hat{w}, \hat{z}$ are the respective solutions of \eqref{wcas} and \eqref{zcas}. We define as before $w(t)=\hat{w}(T-t)$, $z(t)=\hat{z}(T-t)$, $\hat{y}_2(t)=y_2(T-t)$ for
$t \in [0,T]$. Moreover we introduce the solution $\hat{y}_1$ of
$$
\begin{cases}
\hat{y}_{1,tt} -\Delta \hat{y}_1 +c \hat{y}_2=0 \quad \mbox{ in } (0,T) \times \Omega \,,\\
 \hat{y}_1=0 \quad \mbox{ in } (0,T) \times \Gamma \,,\\
(\hat{y}_1,\hat{y}_{1,t})_{| t=0}=(0,0) \quad \mbox{in } \Omega \,.
\end{cases}
$$
\noindent Then, $(\hat{y}_1,\hat{y}_2, \hat{y}_1^{\prime},\hat{y}_2^{\prime})$ is the solution of a cascade system of the form \eqref{cascade1} with
control $v(T-\cdot)$, source term $\xi(T-\cdot)$ and initial data $(\hat{y}_1,\hat{y}_2, \hat{y}_{1,t},\hat{y}_{2,t})_{| t=0}=(0,y_2(T),0,-y_{2,t}(T))$. By definition
of the solutions by transposition of this cascade system (see \eqref{TRANSy1}), we have
\begin{multline*}
0=-\int_0^T\int_{\Omega} c\hat{y}_2 w\,dx\,dt=- \int_0^T\langle \hat{y}_2 \,, cw\rangle_{H_1,H_{-1}}\,dt=
\langle \hat{y}_{1,t}(T), w(T)\rangle_{H_1,H_{-1}} -\\
\langle \hat{y}_1(T), w_t(T)\rangle_{H_2,H_{-2}} -
\Big(
\langle \hat{y}_{1,t}(0), w(0)\rangle_{H_1,H_{-1}} -
\langle \hat{y}_1(0), w_t(0)\rangle_{H_2,H_{-2}}\Big)\,.
\end{multline*}
\noindent Since the initial data for $\hat{y}^1$ are vanishing and since $w(T)=z^0$ and $w_{t}(T)=0$, we have
$$
0=- \int_0^T\langle \hat{y}_2 \,, cw\rangle_{H_1,H_{-1}}\,dt=
\langle \hat{y}_{1,t}(T), z^0\rangle_{H_1,H_{-1}} \quad \forall \ z^0 \in H^1_0(\Omega)\,.
$$ 
\noindent Hence by reflexivity of $H_1$, and density of $H_1=H^1_0(\Omega)$ in $H_{-1}$, we deduce that $\hat{y}_{1,t}(T)=0$.
In a similar way replacing $w$ by $z$ and since $z(T)=0$ and $z_t(T)=-z^1$, we have 
$$
0=- \int_0^T\langle \hat{y}_2 \,, cz\rangle_{H_1,H_{-1}}\,dt=
\langle \hat{y}_{1}(T), z^1\rangle_{H_2,H_{-2}} \quad \forall \ z^1 \in L^2(\Omega)\,.
$$
\noindent Hence by reflexivity of $H_2$ and density of $H=L^2(\Omega)$ in $H_{-2}$, we deduce that $\hat{y}_{1}(T)=0$. We set
$y_1(t)=\hat{y}_1(T-t)$ for $t \in [0,T]$. Then $(y_1,y_2, y_1^{\prime},y_2^{\prime})$ is the solution of  \eqref{cascade1} and satisfies \eqref{insen0}.

The proof of the above results is similar in the boundary control case using the definition of transposition solutions and is left to the reader.
\end{Proof}
\begin{Proof} of Theorem~\ref{thminsens}.
We first consider the case of distributed control. Let $Y^0=(0,y^0,0,y^1) \in X_{-1}^{\ast}$. We consider the bilinear form $\Lambda$ on $X_{-1} $ defined by 
$$
\Lambda(W^T,\widetilde{W}^T)=\int_0^T\int_{\Omega} bw_2 b \widetilde{w_2}\,dx\,dt \,.
$$
and the continuous linear form on $X_{-1}$ defined for all $W^T  \in X_{-1}$
by
\begin{equation}\label{Linsen}
\mathcal{L}(W^T)=
\langle y^1, w_2(0)\rangle_{H,H} - \langle y^0, w_{2,t}(0)\rangle_{H_1,H_{-1}} \,,
\end{equation}
\noindent where $W=(w_1,w_2,w_{1,t},w_{2,t})$ is the solution of 
\begin{equation}\label{Winsen}
\begin{cases}
w_{1,tt} - \Delta w_1=0 \,, \\
w_{2,tt}  - \Delta w_2 + c w_1=0 \,, \\
W_{|t=T}=W^T \,,
\end{cases}
\end{equation}
\noindent and
$\widetilde{W}=(\widetilde{w_1},\widetilde{w_2},\widetilde{w_{1,t}},\widetilde{w_{2,t}})$ is the solution
of \eqref{Winsen} where $W^T$ is replaced by $\widetilde{W}^T$.
We also define a linear form on $X_{-1}$ by
\begin{equation}\label{formxi}
J(W^T)=\int_0^T\int_{\Omega} \xi w_2 \,dx\,dt  \,.
\end{equation}
\noindent Thanks to the inequalities \eqref{inte1U2m} and \eqref{directweak1bis} applied to $Z_2$ where $Z=\mathcal{A}_2^{-1}W$, we deduce that $J$ is continuous on $X_{-1}$.
Thanks respectively to the admissibility inequality\eqref{directweak1bis} and to the observability inequalities \eqref{directweak2bis}-\eqref{directweak3bis}, $\Lambda$ is continuous and coercive on $X_{-1}$ for $T > T_3$. Hence, 
thanks to Lax-Milgram Lemma, there exists a unique $W^T \in X_{-1}$ such that
\begin{equation}\label{HUMinsen}
\Lambda(W^T,\widetilde{W}^T)=-\mathcal{L}(\widetilde{W}^T) - J(\widetilde{W}^T)
\,, \quad \forall \widetilde{W}^T \in X_{-1} \,.
\end{equation}
\noindent We set $v=b w_2$. Then $v \in L^2([0,T];L^2(\Omega))$ and we have by definition of
the solution of \eqref{cascade1} by transposition
\begin{multline*}
\int_0^T \int_{\Omega} v b\widetilde{w_2}\,dx \,dt=
\langle y_{1,t}(T), \widetilde{w_1}(T)\rangle_{H_1,H_{-1}} -
\langle y_1(T), \widetilde{w_{1,t}}(T)\rangle_{H_2,H_{-2}} +
\langle y_{2,t}(T), \widetilde{w_2}(T)\rangle_{H,H} - \\
\langle y_2(T), \widetilde{w_{2,t}}(T)\rangle_{H_1,H_{-1}} - 
\mathcal{L}(\widetilde{W}^T) - J(\widetilde{W}^T)\,, \forall \ 
\widetilde{W}^T \in X_{-1}\,.
\end{multline*}
\noindent On the other hand, we have
$$
\int_0^T \int_{\Omega} v \,b\widetilde{w_2}\,dx\,dt=\Lambda(W^T,\widetilde{W}^T)=-\mathcal{L}(\widetilde{W}^T) - J(\widetilde{W}^T)\,,
$$
\noindent so that,we deduce from these two relations that $Y(T)=(y_1,y_2,y_1^{\prime},y_2^{\prime})(T)=0$. In particular we have $y_1(T)=y_{1,t}(T)=0$ where $(y_1,y_2)$ is the solution of
\eqref{cascade1}. Applying Proposition~\ref{insen4}, we conclude that $v$ insensitizes $\Phi$ along the solutions of \eqref{insen1}. Since
$y_2(T)=y_{2,t}(T)=0$ holds, $v$ is also an exact control for the scalar wave equation \eqref{cascadesing1} solved by $y_2$.

\medskip

We now turn to the case of boundary control. Let $Y_0=(0,y^0,0,y^1)\in X_{1}^{\ast}$.  Thanks to the hidden regularity property induced by the admissibility inequality \eqref{directweak1}, the bilinear form
\begin{equation}\label{Lambdaunbinsen}
\Lambda(U^T,\widetilde{U}^T)= \int_0^T\int_{\Gamma} b \frac{\partial u_2}{
\partial \nu}b \frac{\partial \widetilde{u_2}}{
\partial \nu}\,d\sigma \,dt \,,
\forall \ U^T, \widetilde{U}^T \in X_1 \,,
\end{equation}
\noindent is well-defined and continuous on $X_1$.
We consider the linear form on $X_1$ defined by
\begin{equation}\label{Lunbinsen}
\mathcal{L}(U^T)= 
\langle y^1, u_2(0)\rangle_{H_{-1},H_1} - \langle y^0, u_{2,t}(0)\rangle_{H,H} \,,
\forall \  U^T  \in X_1 \,.
\end{equation}
Thanks respectively to the admissibility inequality \eqref{admissineq} and to the observability inequality
\eqref{obsNSH}, $\Lambda$ is continuous and coercive on $X_1$ for $T > T_3$. Thanks to
\eqref{inte1U2m} and to \eqref{admissineq}, $J$ is continuous on $X_1$. Hence, 
thanks to Lax-Milgram Lemma, there exists a unique $U^T \in X_{1}$ such that

\begin{equation}\label{HUMin}
\Lambda(U^T,\widetilde{U}^T)=-\mathcal{L}(\widetilde{U}^T)- J(\widetilde{U}^T) \,, \quad \forall \ \widetilde{U}^T \in X_1 \,.
\end{equation}
\noindent We set $v_2=b \frac{\partial u_2}{\partial \nu}$. We deduce as for the previous case that $Y(T)=0$. Hence $y_1(T)=y_{1,t}(T)=0$ where $(y_1,y_2)$ is the solution of
\eqref{cascade2}. Applying Proposition~\ref{insen4}, we conclude that $v$ insensitizes $\Phi$ along the solutions of \eqref{insen2}. Since
$y_2(T)=y_{2,t}(T)=0$ holds, $v$ is also an exact control for the scalar wave equation \eqref{cascadesing2} solved by $y_2$.
\end{Proof}

\section{Conclusion and perspectives}
This paper presents several results that solve partially or completely, open questions on the control
of $2$-coupled cascade wave, heat or Schr\"odinger systems. In particular, we give a necessary and sufficient condition on the coupling and control coefficients for the observation of the full cascade hyperbolic system by a single, either bounded or unbounded observation. This leads to results for the existence of insensitizing controls for the scalar wave equation. Such controls are built in a such a way that a given localized measure of the solution, is robust to small unknown variations of the initial data.
We give a complete answer, under the form of a necessary and sufficient condition for the existence of such insensitizing controls. In the case of spatial domains such that their boundaries have no contact of infinite order with its tangent (or are analytic), this necessary and sufficient condition says that the (either internal or boundary) control region and the localized measure region should satisfy both the Geometric Control Condition~\cite{blr92}.

 Our results are based on an extension of the two-level energy method introduced in~\cite{alacras01, sicon03} for $2$-symmetric coupled systems (see also~\cite{alaleaucont, alaleau11} in the case of partially coercive coupling operators), to $2$-coupled cascade systems. Indeed $2$-coupled symmetric systems are conservative, whereas the total energy of $2$-coupled cascade system is not preserved through time. 
 
 The method presented in this paper does not provide the minimal control time for the full coupled system, on the other hand it is a constructive method. 
In~\cite{these-leautaud}, the authors characterize the minimal control time in the case of locally distributed observation in a compact manifold without boundary, by contradiction arguments.
It would be interesting to have characterizations of the minimal control time in the case of boundary observation, and to give its expression in simple geometric conditions. It would also be interesting to characterize this minimal control time by direct and constructive arguments. 

Many other questions are challenging.  For instance, it is important to extend the present results to $n$-coupled cascade systems of coupled PDE's. We have several recent results in this direction but for length reasons and for the sake of clarity, we present these further results on $n$-coupled cascade systems in a different paper \cite{alamultilevel}. These results have for instance applications on the simultaneous control of devices coupled in parallel.

Further generalizations are  to determine whether if it would be possible to extend our positive results to cascade systems with different operators
$A_1$ and $A_2$ for the coupled system \eqref{NSH}. The first results in this direction are given for symmetric coupled systems in
\cite{sicon03} in situations for which $A_1=-r_1 \Delta$, $A_2=-r_2\Delta$  with $r_1 \neq r_2$ (with restrictions on $r_1/r_2$) and under geometric restrictions and for specific forms of couplings. We also provide positive examples with operators of different order such as $A_1=-\Delta$
and $A_2=\Delta ^2$ allowing to couple different PDE's such as wave and plate equations. Negative results for $2$-coupled cascade systems have been obtained in this direction in \cite{these-leautaud} in the torus for wave and heat  one-dimensional cascade systems for two different diffusion coefficients (see also \cite{DLRL}).
Such questions have also been considered for $2$-coupled parabolic cascade systems in~\cite{FCGBDeT10} with constant couplings
and in the case of boundary control. In this paper, the authors  exhibit
a one-dimensional parabolic cascade system with two different operators of the form $A_1= -d_1 \Delta$ and $-d_2 \Delta$, with $d_1\neq d_2$
and $\sqrt{\frac{d_1}{d_2}} \in \mathbb{Q}$, for which null controllability does not hold. They also provide examples for which approximate controllability holds whereas null controllability does not hold. 

One can also consider the same principal operator $A$ in both equations
of \eqref{NSH}, perturbed in each equation by bounded terms of the form $C_{11}u_1$ and $C_{22}u_2$  respectively in the first and second equations, that is
$$
\begin{cases}
u_1^{\prime\prime} + A u_1 +C_{11}u_1= 0 \,,\\
u_2^{\prime\prime} + A u_2+ C_{22}u_2 + C_{21}u_1  = 0 \,,\\
(u_i,u_i^{\prime})(0)=(u_i^0,u_i^1) \mbox{ for }
i=1,2 \,,
\end{cases}
$$
\noindent where $C_{11}$ and $C_{22}$ are bounded operators in the pivot space $H$. In this case the two operators $A_1$ and $A_2$ differ only by bounded perturbations, it would be also interesting to determine necessary and sufficient conditions for observability by a single observation on the second component $u_2$. In particular our results can be generalized to the case $C_{11}=C_{11}^{\ast}=C_{22}$ with $A+C_{11}$ satisfying the assumption $(A1)$ (it is then sufficient to replace $A$ by $A+C_{11}$ in the proof). Extensions to more general coupling terms is an open question.

Several applications, in particular for reaction-diffusion systems, involve nonlinear control systems. In the spirit of~\cite{CGR10}, it would be of great interest to extend the present results to nonlinear coupled systems.

\medskip

\noindent {\bf Acknowledgements}

\medskip

 I would like to thank Luc Miller for helpful discussions on the sharpness of the Geometric Control Condition and on the transmutation method.
I am also very grateful to the referees for their careful reading, and valuable comments and suggestions.


\begin{thebibliography}{abc99xys}
\bibitem{alacras01} {\sc F. Alabau-Boussouira, }{\em Indirect boundary observability of a weakly coupled wave system}. C. R. Acad. Sci. Paris, t. 333, S\'erie I, 645--650 (2001).

\bibitem{sicon03} {\sc F. Alabau-Boussouira, }{\em A two-level energy method for indirect boundary observability and controllability of weakly coupled hyperbolic systems}. Siam J. Control Opt. 42, 871--906 (2003).

\bibitem{alaleaucont} {\sc F. Alabau-Boussouira and M. L\'eautaud, }{\em Indirect controllability of locally coupled systems under geometric conditions}. C. R. Acad. Sci. Paris, t. 349, S\'erie I, 395--400 (2011).

\bibitem{alaleau11} {\sc F. Alabau-Boussouira and M. L\'eautaud, }{\em Indirect controllability of locally coupled wave-type systems and applications}.  Journal de Math\'ematiques Pures et Appliqu\'ees 99, 544--576 (2013).

\bibitem{arxiv} {\sc F. Alabau-Boussouira, }{\em Controllability of cascade coupled systems of multi-dimensional
evolution PDE's by a reduced number of controls}. C. R. Acad. Sci. Paris, t. 350, S\'erie I, 577--582 (2012) (see also arXiv:1109.6645v2).
 
\bibitem{alaamm11}{\sc F. Alabau-Boussouira and K. Ammari, }{\em Sharp energy estimates for nonlinearly locally damped PDEs via observability for the associated undamped system}. J. of Functional Analysis 260, 2424--2450 (2011).

\bibitem{alamultilevel}{\sc F. Alabau-Boussouira, }{\em A hierarchic multi-levels energy method for the control of bi-diagonal and mixed n-coupled cascade systems of PDEÕs by a reduced number of controls}. Submitted (2012).

\bibitem{AKBD06}{\sc F. Ammar-Khodja and A. Benabdallah and C. Dupaix, }{\em Null controllability of some reaction-diffusion systems with one control force}. J. Math. Anal. Appl. 320, 928--943 (2006).

\bibitem{AKBDGB09}{\sc F. Ammar-Khodja and A. Benabdallah and C. Dupaix and M. Gonz\'alez-Burgos, }{\em A {K}alman rank condition for the localized distributed controllability of a class of linear parabolic systems}. J. Evol. Equations 9, 267--291 (2009).

\bibitem{AKBGBT11}{\sc F. Ammar-Khodja, A. Benabdallah, M. Gonz\'alez-Burgos, L. de Teresa, }{\em Recent results on the controllability of linear coupled parabolic problems: a survey}. Mathematical Control and Related Fields 1, 267--306 (2011).

\bibitem{blr92}{\sc C. Bardos, G. Lebeau, J. Rauch, } {\em Sharp sufficient conditions for the observation, control, and stabilization of waves from the boundary}. SIAM J. Control Opt. 30, 1024--1065 (1992).

\bibitem{bodart-fabre95}{\sc O. Bodart and C. Fabre, }{\em Controls insensitizing the norm of the solution of a semilinear heat equation}. J. Math. Anal. and App. 195, 658--683 (1995).

\bibitem{BGBPG04CPDE}{\sc O. Bodart, M. Gonz\'alez-Burgos and R. P\'erez-Garc\'{i}a, }{\em Existence of insensitizing controls for a semilinear heat equation with a superlinear nonlinearity}. Comm. Partial Differential Equations 29 (2004), 1017--1050 (2004).

\bibitem{BGBPG04SICON} {\sc O. Bodart, M. Gonz\'alez-Burgos and R. P\'erez-Garc\'{i}a, }{\em A local result on insensitizing controls for a semilinear heat equation with nonlinear boundary Fourier conditions}. SIAM J. Control Optim. 43, 95--969 (2004).

\bibitem{burq} {\sc N. Burq, } {\em Contr\^olabilit\'e exacte des ondes dans des ouverts peu r\'eguliers}. Asymptot. Anal. 14, 157--191 (1997).

\bibitem{burqgerard} {\sc N. Burq and P. G\'erard, } {\em Condition n\'ecessaire et suffisante pour la contr\^olabilit\'e exacte des ondes}. C. R. Acad. Sci. Paris S\'er I Math 325, 749--752 (1997).

\bibitem{coronbk07}{\sc J.-M. Coron, }{\em Control and nonlinearity}. Mathematical Surveys and Monographs, 136. American Mathematical Society, Providence, RI, (2007).

\bibitem{CGR10}{\sc J.-M. Coron, S. Guerrero, L. Rosier, }{\em Null controllability of a parabolic system with a cubic coupling term}. SIAM J. Control Optim. 48, 5629--5653 (2010).

\bibitem{Dager06}{\sc R. D{\'a}ger, }{\em Insensitizing controls for the 1-{D} wave equation}. Siam J. Control Opt. 45, 1758--1768 (2006).

\bibitem{DLRL}{\sc B. Dehman, J. Le Rousseau and M. L\'eautaud, }{\em Controllability of two coupled wave equations on a compact manifold}. Preprint HAL hal-00686967 (2012).

\bibitem{EZ}{\sc S. Ervedoza and E. Zuazua, }{\em Observability of heat processes by transmutation without geometric restrictions}. Math. Control Relat. Fields 1, 177-187 (2011).

\bibitem{EZ2}{\sc S. Ervedoza and E. Zuazua, }{\em Sharp observability estimates for heat equations}. Arch. Ration. Mech. Anal. 202, 975--1017 (2011).

\bibitem{FCGBDeT10}{\sc E. Fern{\'a}ndez-Cara and M. Gonz{\'a}lez-Burgos and L. de Teresa, }{\em Boundary controllability of parabolic coupled equations}. J. Functional Anal. 259, 1720--1758 (2010). 

\bibitem{GBdT10} {\sc M. Gonz\'{a}lez-Burgos and L. de Teresa, }
{\em Controllability results for cascade systems of {$m$} coupled parabolic {PDE}s by one control force}. Port. Math. 67, 91--113 (2010).


\bibitem{guerrero}{\sc S. Guerrero, }{\em Controllability of systems of Stokes equations with one control force: existence of insensitizing controls}. Annales de l'Institut Henri Poincar\'e Analyse Non Lin\'eaire 24, 1029--1054 (2007).

\bibitem{gueyecras}{\sc M. Gueye, }{\em Uniqueness results for Stokes cascade systems and application to insensitizing controls}. C. R. Acad. Sci. Paris, t. 350, S\'erie I, 831--835 (2012).

\bibitem{gueyeNS}{\sc M. Gueye, }{\em Insensitizing controls for the Navier-Stokes equations}.  Accept\'e pour publication aux Annales IHP.  ArXiv:1207.3255 (2012).

\bibitem{jaffard}{\sc S. Jaffard, }{\em Contr\^ole interne exact des vibrations d'une plaque rectangulaire}.  Port. Math. 47, 423--429 (1990).

\bibitem{DeTK10}{\sc O. Kavian and L. de Teresa, } {\em Unique continuation principle for systems of parabolic equations}. ESAIM COCV 16, 247--274 (2010).

\bibitem{komornik1997}{\sc V. Komornik, }{\em Rapid boundary stabilization of linear distributed systems}. Siam J. Control Opt. 35, 1591--1613 (1997).

\bibitem{lasieckabk02}{\sc I. Lasiecka, }{\em Mathematical control theory of coupled PDEs}. CBMS-NSF Regional Conference Series in Applied Mathematics, 75. Society for Industrial and Applied Mathematics (SIAM), Philadelphia, PA, (2002).

\bibitem{leautaud}{\sc M. L\'eautaud, }{\em Spectral inequalities for non-selfadjoint elliptic operators and application to the null-controllability of parabolic systems}. J. Funct.
Anal. 258, 2739--2778 (2010).

\bibitem{these-leautaud}{\sc M. L\'eautaud, }{\em Quelques probl\`emes de contr\^ole d'\'equations aux d\'eriv\'ees partielles : in\'egalit\'es spectrales, syst\`emes coupl\'es et limites singuli\`eres}. Th\`ese de doctorat de l'universit\'e de Paris 6 (2011).

\bibitem{lions} {\sc J. L. Lions, }{\em Contr\^olabilit\'e exacte et stabilisation  de syst\`emes distribu\'es}. Vol. 1--2, Masson, Paris (1988).

\bibitem{lions89}{\sc J.L. Lions, }{\em Remarques pr\'eliminaires sur le contr\^ole des syst\`emes \`a donn\'ees incompl\`etes}. Actas del Congreso de Ecuaciones Diferenciales y Aplicaciones (CEDYA), Universidad de M\'alaga, 43--54 (1989).

\bibitem{mauffrey}{\sc K. Mauffrey, }{\em On the null controllability of a 3 $x$ 3 parabolic system with non constant coefficients by one or two control forces}. Journal de Math\'ematiques Pures et Appliqu\'ees 99, 187--210 (2013).

\bibitem{miller04}{\sc L. Miller, }{\em Geometric bounds on the growth rate of the null-controllability cost for the heat equation in small time}. J. of Differential Eq. 204, 202--226 (2004).

\bibitem{miller05}{\sc L. Miller, }{\em Controllability cost of conservative systems: resolvent condition and transmutation}. J. of Functional Anal. 218, 425--444 (2005).

\bibitem{miller}{\sc L. Miller, }{\em The control transmutation method and the cost of fast controls}. Siam J. Cont. Opt. 45, 762--772 (2006).

\bibitem{olive}{\sc G. Olive, }{\em Null-controllability for some linear parabolic systems with
controls acting on different parts of the domain and its boundary}. Mathematics of Control, Signals and Systems 23, 257--280 (2012).

\bibitem{phung} {\sc K.-D. Phung, }{\em Observability and control of {S}chr\"odinger equations}. Siam J. Cont. Opt. 40, 211--230 (2001).

\bibitem{RdT11}{\sc L. Rosier and L. de Teresa, }{\em Exact controllability of a cascade system of conservative equations}. C. R. Acad. Sci. Paris, Ser. I 349, 291--296 (2011).

\bibitem{seidman}{\sc T. Seidman, }{\em Two results on exact boundary control of parabolic equations}. Appl. Math. Optim. 11, 145--152 (1984).

\bibitem{tebou08}{\sc L. Tebou, }{\em Locally distributed desensitizing controls for the wave equation}. C. R. Acad. Sci. Paris, t. 346, S\'erie I, 407--412 (2008).

\bibitem{tebou2011}{\sc L. Tebou, }{\em Some results on the controllability of coupled semilinear wave equations: the desensitizing control case}. Siam J. Control Opt. 49, 1221--1238 (2011).

\bibitem{tenentucs}{\sc G. Tenenbaum, M. Tucsnak, }{\em Fast and strongly localized observation for the schr\"odinger equation}. Trans. of the A.M.S. 351, 951--977 (2009).

\bibitem{DeT00} {\sc L. de Teresa, } {\em Insensitizing controls for a semilinear heat equation}. CPDE 25, 39--72 (2000).

\bibitem{DeTZ} {\sc L. de Teresa and E. Zuazua, } {\em Identification of the class of initial data for the insensitizing control of the heat equation}. CPAA 8, 457--471 (2009).
\end{thebibliography}
\end{document}